\documentclass[12pt, twoside, a4paper]{article}
\usepackage[english]{babel}  
\usepackage{float}
\usepackage[utf8]{inputenc}  
\usepackage[T2A]{fontenc}      
\usepackage{amsmath, amssymb, amsthm, amsfonts, amsthm, dsfont,mathrsfs,wrapfig,graphicx}
\usepackage[]{hyperref}
\usepackage{epigraph}
\usepackage{accents}

\usepackage{abstract}
\DeclareMathOperator{\supp}{supp}

\DeclareMathOperator{\pr}{pr}
\DeclareMathOperator{\const}{const}

\DeclareMathOperator{\loc}{{loc}}
\DeclareMathOperator{\Ree}{Re}
\DeclareMathOperator{\Imm}{Im}
\DeclareMathOperator{\Exc}{Asc}

\DeclareMathOperator{\Cyl}{{Cyl}}

\newcommand{\Sp}{\mathcal F}

\DeclareMathOperator{\SL}{SL}
\DeclareMathOperator{\arctg}{arctg}
\DeclareMathOperator{\Maas}{Asc}
\DeclareMathOperator{\tg}{tg}
\DeclareMathOperator{\Op}{Op}
\DeclareMathOperator{\SSS}{S}
\DeclareMathOperator{\spec}{spec}
\DeclareMathOperator{\Meas}{Meas}
\DeclareMathOperator{\Sc}{Sc}
\DeclareMathOperator{\Isom}{Isom}

\newcommand{\eps}{\varepsilon}
\newcommand{\R}{{\mathbb R}}
\newcommand{\HH}{{\mathbb H}}
\newcommand{\emm}{\eta_{\max}}

\newtheorem{theorem}{Theorem}[section]
\newtheorem{lemma}[theorem]{Lemma}
\newtheorem{define}[theorem]{Definition}
\newtheorem{predl}[theorem]{Proposition}
\newtheorem{sled}[theorem]{Corollary}
\newtheorem{informal_statement}[theorem]{Informal Proposition}

\newcommand{\tk}{B}
\newcommand{\tm}{{\tilde m}}
\newcommand{\B}{\mathsf B}
\newcommand{\A}{\mathsf A}

\usepackage{enumitem}

\usepackage{fancyhdr}
\renewenvironment{abstract}
{\small
	\begin{center}
		\bfseries \abstractname\vspace{-.5em}\vspace{0pt}
	\end{center}
	\list{}{%
		\setlength{\leftmargin}{2cm}
		\setlength{\rightmargin}{\leftmargin}%
	}%
	\item\relax}
{\endlist}

\makeatletter
\def\blindfootnote{\gdef\@thefnmark{}\@footnotetext}
\makeatother

\usepackage{mathtools}

\begin{document}

\LARGE 
\Large
\normalsize

\thispagestyle{empty}
\title{\vspace{-1cm}Infinite ascension limit: horocyclic chaos}
\author{Mikhail Dubashinskiy} 
\date{\today}
{\vspace{-6cm}\maketitle}

\blindfootnote{Chebyshev Laboratory, St.~Petersburg State University, 14th Line 29b, Vasilyevsky Island, Saint~Petersburg 199178, Russia.}
\blindfootnote{\hspace{0.2mm}e-mail: \href{mailto://mikhail.dubashinskiy@gmail.com}{\texttt{mikhail.dubashinskiy@gmail.com}}}
\blindfootnote{Research is supported by the Russian Science Foundation grant 19-71-30002.}

\blindfootnote{{Keywords:} \emph{quantum unique ergodicity, raising and lowering operators, singularity propagation}.}

\blindfootnote{JGP keywords: \emph{ergodic theory, quantum dynamical and integrable systems,  Lagrangian and Hamiltonian mechanics.}}

\blindfootnote{\hspace{0.0mm}{MSC 2010 Primary}: 58J51; Secondary:  11F37, 37D40, 74J20.
}

\blindfootnote{DOI: \url{https://doi.org/10.1016/j.geomphys.2020.104053}}

\blindfootnote{This manuscript version is made available under the \href{http://creativecommons.org/licenses/bync-nd/4.0/}{CC-BY-NC-ND 4.0 license}.}


\thispagestyle{empty}

\renewcommand{\abstractname}{}

\begin{abstract}\vspace{-1.5cm}\footnotesize
\noindent {\bf Abstract.} 
What will be if, given a pure stationary state on a compact hyperbolic surface, we start applying \emph{raising} operator every $\hbar$ "adiabatic"{} second? It turns that during adiabatic time comparable to $1$ wavefunction will change as a wave traveling with a finite speed (with respect to the adiabatic time), whereas the semiclassical measure of the system will undergo  a controllable transformation possessing a simple geometric description. If adiabatic time goes to infinity then, by quantized Furstenberg Theorem, the system will become quantum uniquely ergodic.

\vspace{-1.5mm}\hspace{0.5cm}Thus, infinite ascension of a closed system leads to quantum chaos.
\end{abstract}

{\footnotesize
	\tableofcontents
}

\pagestyle{fancy}

\section{Introduction}

Consider a \emph{hyperbolic surface} $X$, that is, a Riemannian manifold of real dimension $2$ having constant Gaussian curvature $-1$, a permanent saddle. We always assume that $X$ is compact and has no boundary. Let $\Delta_X$ be hyperbolic Laplace--Beltrami operator on~$X$. It has purely discrete spectrum due to compactness of $X$. So let $u=u^{(s)}\colon X\to\R$ be Laplace--Beltrami eigenfunctions with $-\Delta_X u^{(s)}=s^2 u^{(s)}$ and $\int_X |u^{(s)}|^2\,d\mathcal A=1$, where $\mathcal A$ denotes the hyperbolic area measure  on $X$ and $s\ge 0$ ranges  discrete set $\sqrt{\spec(-\Delta_X)}$ accumulating to $+\infty$ (hereafter $\spec(-\Delta_X)$ denotes spectrum of operator $-\Delta_{X}$). We mostly drop the superscript~${}^{(s)}$ to simplify the notation.

The well-known Quantum Unique Ergodicity (QUE) conjecture by Rudnick and Sarnak states the uniform distribution in $S^*X$ of the whole sequence 
$\{u^{(s)}_0\}_{s\in\sqrt{\spec(-\Delta_X)}}$. (In particular, this would imply that measures  $(u^{(s)})^2\cdot \mathcal A$ on $X$ converge weak* to normed uniform measure $\mathcal A/\mathcal A(X)$ as $s\to\infty$; the precise meaning of Quantum Unique Ergodicity will be specified in Subsection \ref{subsec:measures} below.) 
This conjecture was formulated in  \cite{RS94}. For arithmetic hyperbolic surfaces, QUE was finally proved in \cite{Lin06}. 

To verify QUE, one needs to show that functions $u^{(s)}$ cannot have \emph{microlocal} singularities in the semiclassical limit.
In this direction, some significant results on deconcentration of eigenfunctions for the general non-arithmetic case were obtained in \cite{An08}, \cite{AnantharamanNonnenmacher} and  \cite{DJ}.

The purpose of our paper is to study the meaning  of Maa\ss {} raising (and lowering) operators in this context; these operators may also be understood  as \emph{creation} (and, respectively, \emph{annihilation}) \emph{operators}. Also, we give quantum counterpart of Furstenberg's Theorem on unique ergodicity of horocyclic flow. The latter assertion turns to be related to Rudnick--Sarnak conjecture by a quantum homotopy given by composition of raisings.

\subsection*{Main results}

We may cover $X$ locally isometrically by  hyperbolic plane $\HH$, the later is  implemented as upper-half plane $\mathbb C^+=\{x+iy\colon y>0, x\in\R\}$. Then, $X$ can be understood as $\Gamma\setminus\HH$ for an appropriate discrete subgroup $\Gamma$ in $\Isom^+(\HH)$, the latter is the group of orientation-preserving isometries of $\HH$; let $F\subset \HH$ be any fundamental domain for $\Gamma$. In what follows, $\tau$ is an integer.

\begin{define}
	\label{def:tau-form}
	A $C^\infty$-function $u\colon \HH\to \mathbb C$ will be called a \emph{$\tau$-form} \emph{(}automorphic with respect to the group $\Gamma$\emph{)} if, for each $z\in\HH$, 
	\begin{equation}
	\label{eq:form_def}
	u(\gamma z) = \left(\frac{cz+d}{c\bar z+d}\right)^\tau  u(z)
	\end{equation}
	for any $\gamma\in\Gamma$ of the form $\gamma(z) = \dfrac{az+b}{cz+d}$,  $z\in\HH\simeq\mathbb C^+$, $a,b,c,d\in\R$. The set of functions $u$ with such automorphy will be denoted by $\Sp^\tau(\Gamma)$. Number $\tau$ is understood as the \emph{degree} of form $u$.  
\end{define}

\noindent For $\tau\in\mathbb Z$, we define \emph{raising operator} $K_\tau\colon C^\infty(\HH)\to C^\infty(\HH)$ by
$$
K_\tau u (z)= 2iy\frac{\partial u}{\partial z} + \tau u(z), ~~ z=x+iy\in\HH, ~~ u\in C^\infty(\HH).
$$
There is also \emph{lowering operator} $L_\tau u (z)= -2iy\dfrac{\partial u}{\partial \bar z} - \tau u(z)$. Since $L_\tau=\bar K_{-\tau}$, the study of lowering operators can be reduced to the study of raisings; we thus restrict ourselves to  considering raising operators in what follows. 

We also deal with  $\tau$-Laplacian 
$D^\tau := -\Delta_\HH+2i\tau y\dfrac{\partial}{\partial x}$. (In coordinates $(x,y)$ in $\HH$, hyperbolic Laplacian takes the form $\Delta_\HH:=y^2\left(\dfrac{\partial^2}{\partial x^2}+\dfrac{\partial^2}{\partial y^2}\right)$.) There is 

\begin{predl}[see \cite{Fay}]
	\label{prop:Maas_properties} 
	Let $\tau$ be an integer, $K_\tau$ and $D^\tau$ be operators defined as above, and $\Gamma$ be an \emph{arbitrary} group of hyperbolic isometries.	
	\begin{enumerate}
		\item Operator $K_\tau$ maps $\Sp^\tau(\Gamma)$ to $\Sp^{\tau+1}(\Gamma)$.
		
		\item Operator $K_\tau$ intertwines $D^\tau$ and $D^{\tau+1}$, that is, 
		$K_\tau D^\tau=D^{\tau+1} K_\tau$.
		
		\item 	If $D^{\tau}u=s^2u$ \emph{(}$u\in C^\infty(\HH)$\emph{)} then $D^{\tau+1}(K_\tau u)=s^2K_\tau u$; in other words, $K_\tau$ takes eigenfunctions of $D^\tau$ to eigenfunctions of $D^{\tau+1}$.
		
		\item If $\Gamma$ is a cocompact group, that is, it has a compact fundamental domain $F$, and $u\in\Sp^\tau(\Gamma)$ is such that  $D^\tau u = s^2 u$ in $\HH$, then 
		$$
		\|K_\tau u\|_{L^2(F)}^2 = \left(s^2+\tau(\tau+1)\right)\cdot\|u\|^2_{L^2(F)}
		$$
		\emph{(}the left-hand side  does not depend on the choice of $F$ because factor $\left(\dfrac{cz+d}{c\bar z+d}\right)^{\tau+1}$ in \emph{(\ref{eq:form_def})} is of unit absolute value, and the same for the expression at the right\emph{)}.
	\end{enumerate}	
\end{predl}

\noindent Now, suppose that a function $u_0=u_0^{(s)}\in \Sp^0(\Gamma)$ is such that $D^0 u_0 = s^2 u_0$, $s\in\R$, that is, $u_0$ is Laplace--Beltrami eigenfunction on $X$; assume that $\|u_0\|_{L^2(X)}=1$. Pick some $\tk\in\R^+$, this parameter is understood as \emph{adiabatic time}. On $C^\infty(\HH)$, define \emph{ascension operator}
$$
\Exc_{0\to \tk}^{(s)} := \frac{K_{[\tk s]-1}}{\sqrt{s^2+([\tk s]-1)\cdot[\tk s]}}\cdot\,\cdots\,\cdot\frac{K_1}{\sqrt{s^2+1\cdot(1+1)}}\cdot\frac{K_0}{\sqrt{s^2+0\cdot(0+1)}}.
$$
Notice that by the explicit form of $K_{\tau}$, numerator in $\Exc_{0\to \tk}^{(s)}$ is something like $\Gamma(K_0+[\tk s])\cdot \Gamma^{-1}(K_0)$, here (and only here) $\Gamma$ is gamma-function, $K_0=2iy\dfrac{\partial}{\partial z}$. 

Now put
$$
u_{\tk} =u_{\tk}^{(s)}:= \Exc_{0\to \tk}^{(s)}\,u_0^{(s)}\in\Sp^{[\tk s]}(\Gamma).
$$ 

\noindent The family $\{\Exc_{0\to \tk}^{(s)}u_0^{(s)}\}_{\tk\ge0}$ will be called \emph{ascension evolution} of function $u_0^{(s)}$, this evolution is parametrized by adiabatic time $\tk$.  Note that application of one raising operator takes about $\hbar=1/s$ adiabatic seconds, this is average length of wave $u_0^{(s)}$. 

By the construction, all the functions $u_\tk$ are of unit norm in $L^2(F)$. Thus, mapping $\Exc_{0\to \tk}^{(s)}$ acting on functions $u\in \mathcal{F}^0(\Gamma)$ with $D^0u=s^2u$ can be understood as an isometric operator, say, with respect to $L^2(F)$-norm (though,  domain of such an operator is usually one-dimensional --- in the case of absence of multiple spectrum). This leads us to the natural desire \emph{to study operator $\Exc_{0\to \tk}^{(s)}$ from the analytical viewpoint}. There is the following 

\begin{informal_statement}
	\label{th:informal_intro}
	Wavefunction $u_\tk = \Exc_{0\to \tk}^{(s)} u_0$ understood as a function of $\tk$ looks like a wave running with bounded speed \emph{(}evaluated with respect to $\tk$\emph{)}.
\end{informal_statement}

\noindent We formalize and verify this observation for cylindrical harmonics via WKB techniques in Subsections \ref{subsec:WKB} and~\ref{subsec:phase_transport}.

And this is also an informal heuristics that if any "elementary"{} wave travels with a bounded speed then we have control on motion of "local frequency spectrum"{} of a function decomposable into such waves; in particular, propagation of singularities can be described. (For quantum Hamiltonian evolutions, this is known as Yu. Egorov Theorem, but our considerations do not fall in this case.) 

The latter observation is formalized by using the notion of \emph{semiclassical measures}. Let $J\subset \sqrt{\spec(-\Delta_{X})}$ be an infinite sequence such that $\{(u_{\tk}^{(s)}, 1/s)\}_{s\in J}$ has a semiclassical measure  $\bar\mu^\tk\in\Meas(T^*\HH)$; to be brief, this notation means that we quantize everything related to function $u_\tk^{(s)}$ using Planck constant $\hbar=1/s$, see Definition \ref{def:semiclassical_measure}. 

For $\hbar=1/s$, 
function $u_\tk$ almost satisfies 
$$\left(-\hbar^2\Delta_{\HH}+2i\tk\hbar y\dfrac{\partial}{\partial x}-1\right) u_\tk=0.$$ 
Roughly speaking, $\Exc_{0\to \tk}^{(s)}$ almost intertwines free-particle quantum Hamiltonian $-\hbar^2\Delta_{\HH}$ and \emph{quantum magnetic Hamiltonian} $-\hbar^2\Delta_{\HH}+2i\tk\hbar y\dfrac{\partial}{\partial x}$. The latter operator is \emph{quantization} of symbol $2H_\tk-\tk^2$, where \emph{classical magnetic Hamiltonian} $H_\tk$ is 
$$
H_\tk := \frac{1}{2}\left((y\xi_1-\tk)^2+(y\xi_2)^2\right)
$$
(see Subsection \ref{subsubsec:quantization} or \cite{Zw} for more details on quantization; $(x,y,\xi_1,\xi_2)$ is the canonical coordinate system in $T^*\HH$, see Subsection \ref{subsec:magnetic}). Therefore, it is natural to consider the identification $\phi_{\tk}\colon T^*\HH\to T\HH$ generated by  $H_\tk$ considered as Hamilton function (see details in Subsection~\ref{subsec:magnetic}). Define $\underline{\mu}^\tk := (\phi_{\tk})_{\sharp}\bar{\mu}^\tk$, push-forward of $\bar{\mu}^\tk$ by~$\phi_{\tk}$.

The very standard fact is that measure $\underline{\mu}^{\tk}$ is supported by the set $S_{\sqrt{\tk^2+1}}\HH$ of tangent vectors of length $\sqrt{\tk^2+1}$. Moreover, $\underline{\mu}^\tk$ is invariant with respect to \emph{$\tk\mbox{-hypercyclic}$ flow} which is motion over curves of constant geodesic curvature $\tk/\sqrt{\tk^2+1}$ passed with constant speed $\sqrt{\tk^2+1}$ (see Subsection~\ref{subsec:geometry}). Also, we have 

\begin{predl}[on invariance]
	\label{prop:invariance}
		\begin{enumerate}[leftmargin=*,parsep=-1pt]
			\item
	Measure $\underline{\mu}^\tk$ is $\Gamma$-invariant measure on $T\HH$, in other words, this is a measure on $TX$. 
	\item As a measure on $TX$, $\underline{\mu}^{\tk}$ does not depend on covering of $X$ by $\HH$ \emph{(}which we have chosen at the beginning of our construction\emph{)}.
\end{enumerate}
\end{predl}

\noindent See the proof at the end of Subsection \ref{subsec:measures}.

Now we give explicit transformation of $TX$ taking $\underline{\mu}^0$ to $\underline{\mu}^\tk$. 
This transformation, naturally, should be a conjugation between geodesic and $\tk$-hypercyclic flows. This is because geodesic and $\tk$-hypercyclic flows are classical versions of quantum free-particle and quantum particle in magnetic field respectively, whereas ascension intertwines dynamics of the latter ones.

More formally, take any $v$ from $TX$ (or from $T\HH$). For $\theta$ real, let $\mathcal R_\theta v$ be rotation of vector $v$ around its basepoint by angle $\theta$ counterclockwise. Also, for $c\ge 0$, let $\Sc_c v$ be vector $v$ scaled $c$ times (with the same basepoint and same direction). We denote by $h_t^0$, $t\in\R$, the geodesic flow acting on $TX$ or on $T\HH$ (see Subsection \ref{subsec:geometry}). Put $$\mathcal T_{\tk}v :=\Sc_{\sqrt{\tk^2+1}}\circ \mathcal{R}_{\pi/2}\circ h^0_{\ln\left(\tk+\sqrt{\tk^2+1}\right)}\circ\mathcal{R}_{-\pi/2} v$$ for $v\in TX$ or $v\in T\HH$. Then $\mathcal T_\tk\colon TX\to TX$ is defined invariantly. 

Our first main result is the following

\begin{theorem}[on finite ascension]
	\label{th:finite_time_shift} Suppose that $J\subset \sqrt{\spec(-\Delta_{X})}$ is such that sequence $\{(u_0^{(s)}, 1/s)\}_{s\in J}$ has semiclassical measure $\bar{\mu}^0$. 	Then $\{(u^{(s)}_\tk,1/s)\}_{s\in J}$ also has  semiclassical measure $\bar{\mu}^\tk$, and  for $\underline{\mu}^0$, $\underline{\mu}^\tk$ defined as above, we have 
	$\underline{\mu}^\tk = \left(\mathcal T_\tk\right)_\sharp \underline{\mu}^0$.
\end{theorem}

\noindent This means that when $\tk$ varies smoothly, $\underline{\mu}^\tk$ also travels smoothly and in a controlled way: initial mass on a vector $v\in TX$ shifts in the direction orthogonal to $v$ to the right of $v$ and by distance $\ln\left(\tk+\sqrt{\tk^2+1}\right)$. 

We may consider mapping $$\tk\mapsto \tk\mbox{-hypercyclic flow}, ~~~ \tk\in\R,$$ as  one-parameter family of classical dynamical systems, or as a homotopy between them. To any such system, there corresponds one-parametric operator group $\exp(tD_{\tk,\hbar}/i\hbar)$, $t\in\R$, $D_{\tk,\hbar}=-\hbar^2\Delta_{\HH}+2i\tk\hbar y\dfrac{\partial}{\partial x}$, acting on $\mathcal F^{\tk/\hbar}(\Gamma)$. Theorem~\ref{th:finite_time_shift} says that ascension operator $\Exc_{0\to \tk}^{(s)}$ not only intertwines quantum Hamiltonians $D_{0,\hbar}$ and $D_{\tk,\hbar}$ but also, in a sense, gives transformation of wave taking stationary states  of $D_{0,\hbar}$ to those of $D_{\tk,\hbar}$; this transformation depends on $\tk$ in a smooth way  --- in the sense of semiclassical measures. Ascension evolution therefore should be understood as \emph{quantum homotopy} between quantum magnetic systems, this homotopy is quantization of homotopy between classical $\tk$-hypercyclic flows on $X$. We have got one more implementation of Bohr principle. Recall that the latter states various kinds of correspondence between classical dynamical systems and their quantizations in the semiclassical limit $\hbar\to 0$.

\medskip

Now, suppose that $0.01$ of mass of some weak* limit of sequence $\{(u_0^{(s)})^2\cdot\mathcal A\}_{s\in J}$ turned to be concentrated on a closed geodesic loop $\gamma$ (to author's best knowledge, such a possibility still is not disproven, at least, it is not prohibited by \cite{An08}, \cite{AnantharamanNonnenmacher} and \cite{DJ}); then $0.01$ of mass of weak* limit of the corresponding sequence $\left\{\left|u_\tk^{(s)}\right|^2\cdot\mathcal A\right\}_{s\in J}$ will be concentrated on the two $\tk$-hypercycles obtained by shifting $\gamma$ by distance $\ln\left(\tk^2+\sqrt{\tk^2+1}\right)$ in the left- and right-normal directions. If $\tk\to+\infty$ then both these $\tk$-hypercycles become long, close to \emph{horocycles} and thence almost uniformly distributed in $X$ and in $S_1X$. This is just by Furstenberg's Theorem   which says that, for compact $X$, there is unique Borel probability measure supported by spherical bundle $S_1 X$ and invariant under the action of horocyclic flow; this property is known as \emph{unique ergodicity} of horocyclic flow (see \cite{Furstenberg73} and also dynamical proof in \cite{Marcus}). 

So, the following question is natural: \emph{what will be if we let $\tk$ in $\Exc_{0\to \tk}^{(s)}$ go to infinity?} This is: what is the behavior of 
$$
U_{s,\tau}=\frac{K_{\tau-1}}{\sqrt{s^2+(\tau-1) \cdot\tau}}\cdot\,\cdots\,\cdot\frac{K_1}{\sqrt{s^2+1\cdot(1+1)}}\cdot\frac{K_0}{\sqrt{s^2+0\cdot(0+1)}}\, u^{(s)}_0
$$
if $\tau\in\mathbb N$ grows faster than $s$? This question can be answered.   We have the following 

\begin{theorem}[infinite ascension limit]
	\label{th:horrorcyclic_chaos}
	For any sequences $$\{s_n\}_{n=1}^\infty\subset\sqrt{\spec\left(-\Delta_X\right)}\setminus\{0\},  ~~  \{\tau_n\}_{n=1}^\infty\subset\mathbb N$$ with $\dfrac{\tau_n}{s_n}\xrightarrow{n\to\infty}\infty$, semiclassical measure of sequence $\{(U_{s_n,\tau_n}, 1/\tau_n)\}_{n=1}^\infty$ is $\dfrac{\left(\phi_{1}^{-1}\right)_\sharp\mu_L}{\mathcal A(X)}$. 
	Here, $\mu_L$ is the uniform Liouville measure on $S_1\HH$.
	
	In other words,  $\{U_{s_n,\tau_n}\}_{n=1}^\infty$ is QUE sequence. Infinite ascension of closed system leads to quantum chaos. 
\end{theorem}

\noindent 
A similar result was obtained earlier by Zelditch in \cite{Ze92}. In that paper, the exposition is given in terms of representation theory. For the convenience of reader, we give a proof in the PDE language in this article. Note also that Zelditch takes the limit over a ladder with fixed eigenvalue whereas, in Theorem \ref{th:horrorcyclic_chaos}, square root of eigenvalue $s$ can go to $\infty$ --- just slower than $\tau$. 

Notice that we use quantization at level $1/\tau_n$ in Theorem \ref{th:horrorcyclic_chaos}, this means that, under assumptions of this Theorem, average wavelength of function $U_{s_n,\tau_n}$ is comparable to $1/\tau_n$. 
This is where we make use of compactness of $X$: in fact, Theorem \ref{th:horrorcyclic_chaos} is a quantum version of Furstenberg Theorem and is derived from the latter. That's because, under assumption $s/\tau\to 0$, function  $U_{s,\tau}$ satisfies $\left(-\hbar^2\Delta_{\HH}+2i y\hbar\dfrac{\partial}{\partial x}\right)U_{s,\tau}=o(1)\cdot U_{s,\tau}$  with $\hbar=1/\tau$, and Hamiltonian dynamics given by the symbol of the operator at the left is horocyclic flow (at the corresponding energy level).

\medskip

This paper arose from an attempt to prove Rudnick--Sarnak QUE conjecture. As it was mentioned above, the family of operators $\{\Exc_{0\to \tk}^{(s)}\}$ with $\tk$ increasing from $0$ to $+\infty$ may be understood as a {quantum homotopy}. 
By Theorem \ref{th:finite_time_shift}, this homotopy preserves chaoticity when adiabatic time $\tk$ ranges a \emph{finite} real interval. Also, this homotopy reaches a certainly chaotic system, the quantization of horocyclic flow (Theorem \ref{th:horrorcyclic_chaos}); it does reach but for an infinite time, this does not provide chaos for the initial system. Thus, this construction does not prove Rudnick--Sarnak conjecture.

\section{Classical and quantum magnetic dynamics}

Until Section \ref{sec:finite_time}, we mostly deal with Lobachevsky hyperbolic plane $\mathbb H$  implemented as $\mathbb C^+$ endowed with the Riemannian metric $ds^2=\dfrac{dx^2+dy^2}{y^2}$; also,  $d\mathcal A=\dfrac{dx\,dy}{y^2}$ is hyperbolic volume element on $\HH$. Further, $T\mathbb H$ and $T^*\mathbb H$ are tangent and cotangent bundles over~$\mathbb H$. For $r\ge0$, let $S_r\mathbb H \subset T\mathbb H$ be the set of all tangent vectors of length $r$.

The very standard \emph{geodesic flow} will be denoted by $h^0_t$, $t\in\R$, this is the one-parameter group acting on $T\mathbb H$. From the physical viewpoint, this is a motion of a free classical particle on $\mathbb H$.

\subsection{Hypercyclic and horocyclic flows as magnetic dynamics}

\label{subsec:geometry}
Now  suppose that our particle has unit charge and mass. Consider uniform magnetic field on $\mathbb H$; its intensity will be denoted by $\tk\in\R$ throughout all the paper. We may think that this field is oriented as "the positive normal field to $\mathbb H$"{}. If $\tk>0$ then the trajectory of the particle in such a field starts curving to the right. Depending on initial speed of particle and on the intensity of magnetic field, the trajectory of the particle can be either a 
geodesic line, a \emph{hypercycle}, a \emph{horocycle} or a circle in hyperbolic metric. In either case, the absolute value of speed of particle remains constant under magnetic dynamics. We are interested in the first three kinds of such curves.

For $\tk\ge 0$, a \emph{$\tk$-hypercycle} is a parametrized curve on $\HH$ of constant geodesic curvature $\dfrac{\tk}{\sqrt{\tk^2+1}}$ curving to the right and passing with the constant speed $\sqrt{\tk^2+1}$. An equivalent definition is as follows: 1. curve $t\mapsto \left(\dfrac{\tk}{\sqrt{\tk^2+1}}\cdot e^t, \dfrac{1}{\sqrt{\tk^2+1}}\cdot e^t\right)$, $t\in\R$, in coordinates $(x,y)$ in $\HH$ is a $\tk$-hypercycle; 2. any shift of this curve by an orientation-preserving isometry of $\HH$ is also a $\tk$-hypercycle. (In $\mathbb C^+$, such curves are Euclidean circles intersecting absolute line $\R$ under angle $\dfrac{\pi}{2}-\arctg\tk$.) Hypercycles are also called \emph{hypercircles} or \emph{equidistant curves}. The latter is because they are equidistant from geodesics. Namely, if $t\mapsto\gamma_0(t)\in\HH$, $t\in\R$, is a geodesic passed with unit speed then parametrized curve $\gamma_\tk(t)$ defined as basepoint  of vector $h^0_{\ln\left(\tk+\sqrt{\tk^2+1}\right)}\circ\mathcal{R}_{-\pi/2}\left(\gamma'_0(t)\right)$ (see Introduction) is a $\tk$-hypercycle. This is easily checked by a direct computation.

A {(}right{)} \emph{horocycle} on Lobachevsky plane $\mathbb H$ is a parametrized curve of constant geodesic curvature $1$ curving to the right and passed with the \emph{unit} speed. An equivalent definition is: 1. the curve $t\mapsto (-t,1)$, $t\in\R$, in $(x,y)$-coordinates in $\HH$ is a horocycle, 2. any shift of this curve by an isometry of $\HH$ is also a horocycle.

Notice that if we reparametrize $\tk$-hypercycles such that they will be passed with unit speed then the obtained curves will tend to horocycles as $\tk\to+\infty$.

For any vector $v$ from $S_{\sqrt{\tk^2+1}}\HH$ (or from $S_1\HH$) there exists a unique $\tk$-hypercycle (respectively, a unique horocycle) parametrized as $t\mapsto \gamma(t)$, $t\in\R$, with $\gamma'(0)=v$. Put $h^\tk_t v := \gamma'(t)$ (respectively, $h^\infty_t v := \gamma'(t)$ for horocycles). In such a way, we have defined \emph{$\tk$-hypercyclic flow} $h^\tk_t\curvearrowright S_{\sqrt{\tk^2+1}}\HH$ and, respectively, \emph{horocyclic flow} $h^\infty_t\curvearrowright S_{{1}}\HH$. These flows are also well-defined on $S_{\sqrt{\tk^2+1}}X$ and $S_1 X$ respectively.  Notice that some authors define horocyclic flow so as basepoint of vector $v$ moves in the direction orthogonal to $v$; that formalization is good for matrix calculations. But we prefer formalism originating in physical intuition.

We have already mentioned that the flow $h_t^\infty$ is uniquely ergodic. 
The flow $h_t^\tk$ is conjugated to the geodesic flow $h^0_t$ (by the mapping $\mathcal T_\tk$ defined in the Introduction); the latter, $h^0_t$, is known to be just ergodic (but is Anosov-type instead and has positive entropy equal to $1$).

\subsection{Magnetic Hamiltonian. Quantization}
\label{subsec:magnetic}
The flows $h_t^\tk$ and $h_t^\infty$ can be defined via Hamiltonian. Let $(x,y,\xi_1,\xi_2)$ be canonical coordinates in $T^*\mathbb H$ where $z=x+iy\in\HH$, $\xi_1$ is conjugate to $x$, $\xi_2$ is conjugate to~$y$. 

The motion  of a classical particle with unit charge and mass in the magnetic field of intensity $\tk$ has Hamiltonian 
$$
H_\tk := \frac{1}{2}\left((y\xi_1-\tk)^2+(y\xi_2)^2\right).
$$
Denote by $\Xi_\tk$ the Hamiltonian vector field given by $H_\tk$, and by  $\exp t\Xi_\tk\colon T^*\HH\to T^*\HH$ ($t\in\R$) the Hamiltonian flow generated by $H_\tk$. 

Any Hamilton function defines a mapping from the cotangent bundle to the tangent bundle which  often turns to be a bijection (see, e.g., \cite{Takhtajan}). Let's check our case. Any vector $v$ tangent to $\HH$ at a point $z=x+iy\in\HH\simeq\mathbb C^+$ can be written as $v=v_x \dfrac{\partial}{\partial x}+v_y \dfrac{\partial}{\partial y}$; so let us take $(x,y,v_x,v_y)$ as a coordinate system in $T\HH$. Length of $v$ is then given by $|v|= y^{-2}\cdot(v_x^2+v_y^2)$.
Define mapping $\phi_\tk\colon T^*\HH\to T\HH$: put 
$$
\phi_\tk(x,y,\xi_1,\xi_2) := \left(x,y,y^2\xi_1-\tk y,y^2\xi_2\right)=\left(x,y, \frac{\partial H_\tk}{\partial \xi_1},\frac{\partial H_\tk}{\partial \xi_2}\right)\in T\HH.
$$
Here, we used the two coordinate systems given just above. The inverse mapping $\phi_\tk^{-1}\colon  T\HH\to T^*\HH$ is
$$
\phi_\tk^{-1}(x,y,v_x,v_y) = \left(x,y,\frac{v_x+\tk y}{y^2},\frac{v_y}{y^2}\right).
$$

\begin{predl}[see also \cite{Sunada}]
	\label{predl:t_tstar_change}
	Let $C\in\R$ be a scalar.
	\begin{enumerate}[leftmargin=*,parsep=-1pt]
		\item Set $\{H_\tk=C\}$ is invariant with respect to the flow $\exp t\Xi_\tk$. 
		\item $\phi_\tk\left(\{H_\tk=C\}\right)=\left\{(z,v)\colon |v|=\sqrt{2C}\right\}$. In particular, $$\phi_\tk\left(\left\{H_\tk=\dfrac{\tk^2+1}{2}\right\}\right)=S_{\sqrt{\tk^2+1}}\HH.$$ 
		\item For any $t\in\R$, $\exp t\Xi_\tk = \phi_\tk^{-1}h^\tk_t \phi_\tk$ on $\{H_\tk=(\tk^2+1)/2\}$.
		\item For any $t\in\R$, $\exp t\Xi_1=\phi_1^{-1}h^\infty_t\phi_1$ on $\{H_1=1/2\}=\phi_1^{-1}(S_1\HH)$.
\end{enumerate}
In other words, flows $h^\tk$ and $h^\infty$ are conjugated to restrictions of $\exp t\Xi_\tk$ \emph{(}respectively, of~$\exp t\Xi_1$\emph{)} to the appropriate level sets of Hamiltonians $H_{\tk}$ and $H_1$ respectively.
\end{predl}

\noindent {\bf Proof.} First claim is standard, second one is obvious. 
Third and fourth are verified by a direct computation. $\blacksquare$

\medskip

\noindent  
Now, instead of hyperbolic plane $\HH$, consider arbitrary hyperbolic surface $X$. The flows $h^\tk_t$ (and $h^\infty_t$) are well-defined on $S_{\sqrt{\tk^2+1}}X$ (and $S_1X$, respectively). Moreover, one can define motion in the uniform magnetic field on $X$ \emph{locally} via appropriate Hamilton function. Nevertheless, such a Hamiltonian cannot be defined globally on a cotangent bundle over \emph{compact} hyperbolic surface with no boundary. This is because the motion in the magnetic field has a plenty of Hamiltonians defined locally on the \emph{cotangent} bundle $T^*X$ (but the flow on $TX$ remains the same because identifications between $TX$ and $T^*X$ given by different Hamilton functions are different). 

From the physical viewpoint, classical magnetic field on $\HH\simeq \mathbb C^+$ is given by a $2$-form $\B=\tk y^{-2}dx\wedge dy$, while Hamiltonian of motion in such a field is  
$$
\frac12\cdot\left|\xi_1-A_1, \xi_2-A_2 \right|^2_{T^*\HH}  = \frac12 \cdot y^2\left((\xi_1-A_1)^2+(\xi_2-A_2)^2\right)
$$
where $\A=A_1 dx+ A_2 dy$ is any primitive of $\B$, that is, $d\A=\B$ (we took $\A=\tk y^{-1}dx$ in order to define our $H_\tk$). Such a primitive can be taken in many ways, up to an exact form; so we cannot expect that we will succeed in taking such a primitive on $X$ as a single-valued form. This difficulty is known as \emph{gauge invariance problem}, see, e.g.,~\cite{Landafshitz}.

Well, suppose that $\A$ is an one-form on $X$ and $d\A=\B$ where the latter $\B$ is $\tk$ times volume form on $X$. But then, by Stokes' Theorem, $\tk\cdot\mathcal A(X)=\int_X \B=\int_X d\A=\int_{\partial X} \A=0$ since $X$ has no boundary; we meet a contradiction if $\tk\neq 0$. Thus, \emph{constant magnetic field is physically impossible on a compact surface}. In other words, we are going to quantize a physically impossible system (more precisely, its energy level). This can be done by replacing wavefunctions by \emph{$\tau$-forms}, tensors with special automorphy.

These tensors have been already defined at Introduction (Definition \ref{def:tau-form}). In some papers this object is called a form of weight $2\tau$.  Note that $\left(\dfrac{cz+d}{c\bar z+d}\right)^\tau$ in Definition~\ref{def:tau-form} does not change if we change signs of $c$ and $d$, that is, if we represent $\gamma$ by  opposite matrix in $\SL(2,\R)$. The assumption $\tau\in\mathbb Z$ simplifies this Definition. But one can get rid of this restriction and consider $\tau$-forms for arbitrary $\tau\in\R$ as in \cite{Fay}. 

Definition  \ref{def:tau-form} is coordinate-dependent, it relies heavily on representation of $z\in \mathbb C^+$ as $z=x+iy$. The same concerns more or less all objects that we define on the cotangent bundle; we will arrive to invariantly defined objects on \emph{tangent} bundle by applying coordinate-dependent identifications $\phi_{\tk}$, see Subsection \ref{subsec:measures} below.

The following operator also depends on coordinate system (or on the covering of hyperbolic surface $X$ by Lobachevsky plane $\HH$). 

\begin{define}
	In $\HH\simeq\mathbb C^+$, define \emph{magnetic Laplacian} \emph{(}or \emph{quantum magnetic Hamiltonian}\emph{)} for magnetic field of uniform intensity $\tk\in\R$ as
$$
D_{\tk,\hbar} := -\hbar^2\Delta_{\HH}+2i\tk\hbar y\dfrac{\partial}{\partial x}.
$$
Here, $\hbar$ is some positive number understood as Planck constant, second-order term $\Delta_{\HH} := y^2\left(\dfrac{\partial^2}{\partial x^2}+\dfrac{\partial^2}{\partial y^2}\right)$ is hyperbolic Laplacian.
\end{define}

\noindent This operator obviously commutes with homotheties $z\mapsto z\cdot e^l$, $z\in\mathbb C^+$, $l\in\R$. Second term in $D_{\tk,\hbar}$, a derivative with respect to the vector field of unit length, arose from the influence of the magnetic field  while the first term is the usual quantum Hamiltonian of a free particle. We see easily that the principal symbol of $D_{\tk,\hbar}$ is $2H_\tk-\tk^2$ up to $O_{\mathcal S^1}(\hbar)$ corrections 
(at least locally; $\mathcal S^1=\mathcal S^1(T^*\HH)$ is the space of Kohn--Nirenberg symbols of the first order, see~\cite{Zw}). Notice also that $\langle D_{\tk,\hbar} u, v\rangle_{L^2(\HH)} = \langle u, D_{\tk,\hbar} v\rangle_{L^2(\HH)}$ if $u,v\in C^\infty(\HH)$ and at least one of these functions is compactly supported.

As we have already mentioned in Introduction, if $\tk/\hbar$ is an integer then $D_{\tk,\hbar}$ acts on the space $\mathcal F^{\tk/\hbar}(\Gamma)$ for any $\Gamma < \Isom^+(\HH)$. Thus, $\mathcal F^{\tk/\hbar}(\Gamma)$ is natural space for quantization of (physically impossible) motion in the magnetic field of intensity $\tk$ on $X=\Gamma\setminus\HH$. Notice that, when $\hbar\to 0$, weight $\tk/\hbar$ goes to infinity; that is, forms from $\mathcal F^{\tk/\hbar}(\Gamma)$ twist faster and faster under change of local conformal coordinates in $X=\Gamma\setminus\HH$. This is crucial for Lemma \ref{lemma:invariance} on invariance of semiclassical measures on tangent space. Now we pass to the definition of these measures.

\subsection{Semiclassical measures and their invariance}

\label{subsec:measures}

The following Definition provides  us with the main tool used to describe the behavior of waves with small wavelengths.

\begin{define}
	\label{def:semiclassical_measure}
	Let $\Omega \subset \HH$ be some open set, functions $v_k$ belong to $L^2(\Omega)$ $(k=1,2,\dots)$, and $\hbar_1, \hbar_2, \dots$ be positive scalars tending to zero. Suppose also that $\sup\limits_{k\in\mathbb N} \|v_k\|_{L^2(\Omega)} < \infty$. We say that a non-negative measure $\mu$ on $T^*\Omega$ is \emph{semiclassical measure} of the sequence $\{(v_k, \hbar_k)\}_{{k\in\mathbb N}}$ if 
\begin{equation}
\label{eq:Wigner_def}
\left\langle(\Op_{\hbar_k} a) v_{k}, v_{k}\right\rangle_{L^2(\Omega)} \xrightarrow[k\to\infty]{} \displaystyle\int_{T^*\Omega} a\,d\mu
\end{equation} for any $a\in C^\infty_0(T^*\Omega)$. Here, linear operator $\Op_{\hbar_k} a\colon L^2(\Omega)\to C_0^\infty(\Omega)$ is the \emph{standard quantization} of symbol $a$ \emph{(}see Subsection \ref{subsubsec:quantization}\emph{)}.
\end{define}

\noindent (If we speak of semiclassical measure of something like $\{(v^{(s)},1/s)\}_{s\in J}$ then it means that limit relation is understood for $s\to\infty$ along $J$.)

Semiclassical measure is sometimes called Wigner measure or Wigner transform of sequence of wavefunctions (don't confuse to Wigner semicircle law!). In Definition \ref{def:semiclassical_measure}, function $a$ on $T^*\Omega$ is understood as \emph{classical observable} used to test the distribution~$\mu$, which, in turn, is understood as a distribution on classical particles. Operator $\Op_\hbar a$ is \emph{quantum observable} applied to wavefunctions; this operator in a sense inherits properties of its symbol $a$ when $\hbar$ is small. Thus, in Definition \ref{def:semiclassical_measure}, wavefunction $v_k$  with small wavelength comparable to $\hbar_k$ gives rise to a distribution on wavevectors. These (co)vectors are local frequencies of $v_k$ scaled $\hbar_k$ times, and they  are identified with classical particles, that is, with points in $T^*\Omega$. See \cite{Zw} for more on semiclassical limits.

A weak*-type argument leads to the following conclusion. \emph{For any sequence of functions  $v_1, v_2, \dots$  bounded uniformly in $L^2(\Omega)$ and for any sequence $\hbar_1, \hbar_2, \dots$ of positive numbers going  to zero, there exists an infinite subsequence of indices $J\subset\mathbb N$ such that the sequence $\{(v_k, \hbar_k)\}_{{k\in J}}$ has a semiclassical measure.} Also, this measure is always non-negative: this is because of almost-positivity of all the reasonable quantization procedures.

Now, we return to the eigenfunctions. Recall that $u_0\in\mathcal F^0(\Gamma)$ is such that $-\Delta_{\HH} u_0 = s^2 u_0$ and $\|u_0\|_{L^2(F)}=1$ where $F$ is any fundamental domain for $\Gamma$; and we put $u_\tk = u_\tk^{(s)}:=\Exc^{(s)}_{0\to \tk} u_0$, operator $\Exc^{(s)}_{0\to \tk}$ being defined in Introduction. Then $u_\tk\in\mathcal F^{[\tk s]}(\Gamma)$ and $\|u_\tk\|_{L^2(F)} = 1$. Therefore, for fixed $\tk\in\R$, we may assume that $J\subset\sqrt{\spec(-\Delta_X)}$ is an infinite  sequence such that both sequences  $\{(u_0^{(s)}, 1/s)\}_{{s\in J}}$ and $\{(u_\tk^{(s)}, 1/s)\}_{{s\in J}}$ have  semiclassical measures which we have already denoted by $\bar\mu^0$ and $\bar\mu^\tk$, respectively.

Rudnick--Sarnak Quantum Unique Ergodicity conjecture mentioned in the Introduction states that $\bar\mu^0$ is the uniform Liouville measure on the set of length $1$ covectors over $\HH$ (or over $X$). To the author's best knowledge, this question is still open. By Lemma~\ref{lemma:Wigner_usual} below and by ergodicity of geodesic flow over $X$, for QUE it is enough to show that $\bar{\mu}^0$ is absolutely continuous with respect to coordinates in the set $\{H_0=1/2\}$.

In this paper, we just study relation between measures $\bar\mu^\tk$ for different $\tk$. By Proposition \ref{prop:Maas_properties}, $D^{[\tk s]} u_\tk =\left(-\Delta_\HH+2i[\tk s] y\dfrac{\partial}{\partial x}\right)u_\tk= s^2 u_\tk$. For convenience, we take $\hbar=\hbar(s) := \dfrac{\tk}{[\tk s]}$ and write
\begin{equation}
\label{eq:u_1_eigen}
D_{\tk,\hbar(s)} u_\tk = \left(-\hbar^2(s)\Delta_{\HH}+2i\tk\hbar(s) y\dfrac{\partial}{\partial x}\right) u_\tk = (1+c(s)) u_\tk, ~~ c(s)\in\R, ~~  c(s)\xrightarrow{s\to \infty}0.
\end{equation}
Note that, by Calderon--Vailliancourt Theorem, $\bar\mu^{\tk}$ is also semiclassical measure for the sequence $\{(u_\tk^{(s)}, \tk/[\tk s])\}_{s\in J}$.

In the proof of the following Lemma and in the rest of the paper, $\langle\cdot, \cdot\rangle$ is $\langle\cdot, \cdot\rangle_{L^2(\HH)}$.

\begin{lemma}[Standard facts on semiclassical measure] 
	\label{lemma:Wigner_usual}
	Let $\bar\mu^\tk$ be semiclassical  measure for the sequence $\{(u^{(s)}_\tk,1/s)\}_{s\in J}$ \emph{(}or $\{(u^{(s)}_\tk,\hbar(s))\}_{s\in J}$\emph{)} satisfying \emph{(\ref{eq:u_1_eigen})}. 
	\begin{enumerate}
		\item Measure $\bar\mu^\tk$ is supported by the set  $\{2H_\tk-\tk^2=1\}=\{H_\tk=(\tk^2+1)/2\}\subset T^*\HH$.
		\item Measure $\bar\mu^\tk$ is invariant with respect to the $\tk$-hypercyclic flow $\exp t\Xi_\tk$ acting on $T^*\HH$ and generated by the Hamiltonian $H_\tk$ \emph{(}and restricted to the energy level $\{H_\tk=(\tk^2+1)/2\}$\emph{)}.
	\end{enumerate}
\end{lemma}

\noindent {\bf Proof.}  First assertion is a direct consequence of \cite[Theorem 5.3]{Zw}. For the second one, we slightly modify the proof of \cite[Theorem 5.4]{Zw}. 
Take any classical real-valued observable $a\in C_0^\infty(T^*\HH)$. Put  $A=\Op^{\R^d}_{\hbar(s)} a$, this is the standard $\R^2\mbox{-quantization}$ in chart $(x,y)$, see (\ref{eq:standard_q_Rd}) in Subsection \ref{subsubsec:quantization} below. Let $\Omega$ be projection of $\supp a$ from $T^*\HH$ to $\HH$. Pick a smooth function $\psi\in C^\infty_0(\HH)$ equal $1$ \emph{near} $\Omega$. We have $D_{\tk,\hbar(s)}(\psi u_\tk) = D_{\tk,\hbar(s)} u_\tk=(1+c(s)) u_\tk$ near $\Omega$. Having in mind symmetricity of~$D_{\tk,\hbar(s)}$, write 
\begin{multline*}
\langle\left[A, \psi D_{\tk,\hbar(s)}\right]\psi u_\tk,u_\tk\rangle=
\langle A \psi D_{\tk,\hbar(s)}\psi u_\tk,u_\tk\rangle-\langle\psi D_{\tk,\hbar(s)}  A\psi u_\tk,u_\tk\rangle=\\=
\langle A \psi D_{\tk,\hbar(s)}\psi u_\tk,u_\tk\rangle-\langle A\psi u_\tk,D_{\tk,\hbar(s)}\psi u_\tk\rangle=
\langle A \psi D_{\tk,\hbar(s)}\psi u_\tk,u_\tk\rangle-(1+c(s))\cdot\langle A\psi u_\tk,u_\tk\rangle=\\=
\langle A \psi \left(D_{\tk,\hbar(s)}\psi -(1+c(s))\right)u_\tk,u_\tk\rangle=O(\hbar^\infty(s)).
\end{multline*}
The latter is because $\psi\left(D_{\tk,\hbar(s)}\psi -(1+c(s))\right)u_\tk=0$ near $\Omega$ whereas $a$ vanishes at covectors with basepoints outside of $\Omega$; we used pseudo-locality of pseudodifferential operators (\cite[p. 211]{Zw}). Since 
\begin{multline*}
[A,\psi D_{\tk,\hbar(s)}] = \dfrac{\hbar(s)}{i}\Op_{\hbar(s)}\{a,\psi(2H_\tk-\tk^2)\}+O_{L^2\to L^2}(\hbar^2(s))=\\=\dfrac{2\hbar(s)}{i}\Op_{\hbar(s)}\{a, H_\tk\}+O_{L^2\to L^2}(\hbar^2(s))
\end{multline*}
($\{\cdot,\cdot\}$ being Poisson brackets), we conclude that $\left\langle\Op_{\hbar(s)}\{a, H_\tk\} \,\psi u_\tk,u_\tk\right\rangle=O(\hbar(s))$ and, by limit pass, that  $\int_{T^*\HH}\{a,H_\tk\}\,d\bar\mu^\tk=0$. This relation holds for any $a\in C^\infty_0(T^*\HH)$, but this implies the invariance of $\bar\mu^\tk$ with respect to~$\exp t\Xi_\tk$. $\blacksquare$

\medskip

\noindent {\bf Remark.} We multiplied $u_\tk$ by cut-off $\psi$ in order to localize quantizations and also to be able to apply $O_{L^2\to L^2}$ estimate to compactly supported $L^2$ function. To be perfect, the same should be done in the following proof of Lemma~\ref{lemma:invariance} below. There, we omit this preparatory step to avoid overcharging the exposition.

\medskip

Our next goal is to prove a kind of invariance of semiclassical measure $\bar\mu^{\tk}$ with respect to $\Gamma$. We have a good chance to succeed since $\tk$-hypercyclic flow is well-defined on $TX$ and this observation has to have some quantum counterpart. 

If $\gamma\colon\HH\to\HH$ is a diffeomorphism then its differential $D\gamma$ is a diffeomorphism of $T\HH$. Thus it makes sense to speak of measures on $T\HH$ invariant with respect to some group of hyperbolic isometries.

\begin{lemma}
	\label{lemma:invariance} 
	Suppose that $\tk\in\R$, $\Gamma<\Isom^+(\HH)$ and that $\hbar$ ranges some set $J\subset \R^+$ accumulating to zero such that $\tk/\hbar$ is always an integer.
	
	Let also $v_\hbar\in\mathcal F^{\tk/\hbar}(\Gamma)$, $\hbar\in J$, be $\Gamma$-automorphic  forms of degrees $\tk/\hbar$ respectively, and suppose that $v_\hbar$ are 
	bounded uniformly with respect to $\hbar$ in $L^2$ on any compact subset of $\HH$.
	
	If $\bar\mu\in\Meas(T^*\HH)$ is  semiclassical measure for the sequence $\{(v_\hbar,\hbar)\}_{\hbar\in J}$ then $\underline{\mu}:=(\phi_{\tk})_{\sharp}{\bar\mu}$ is {$\Gamma$-invariant} measure on $T\HH$.
\end{lemma}

\noindent {\bf Proof.}
We make use of automorphy property of function $v_\hbar$:
\begin{equation}
\label{eq:v_spin}
v_\hbar(\gamma z) = \left(\dfrac{cz+d}{c\bar z+d}\right)^{\tk/\hbar} v_\hbar(z) ~~~ (z\in\mathbb C^+)
\end{equation}
for any $\gamma\in\Gamma$ of the form $\gamma z =\tilde z= \dfrac{az+b}{cz+d}$; we are going to study push-forward of $(\phi_{\tk})_{\sharp}{\bar\mu}$ by mapping $D\gamma$.

Since $\gamma'(z) = \dfrac{1}{(cz+d)^2}$,  mapping $\gamma$  transforms covectors by pull-back as following:
$$
(\tilde{x}+i\tilde y, \tilde{\xi_1}, \tilde\xi_2) \mapsto
(\gamma^{-1}(\tilde{x}+i\tilde y), \alpha\tilde{\xi}_1+\beta\tilde\xi_2, 
-\beta\tilde{\xi}_1+\alpha\tilde\xi_2), ~~~ (\tilde{x},\tilde y, \tilde{\xi_1}, \tilde\xi_2) \in T^*\HH,
$$
where $\alpha=\alpha(z)=\Ree\left(\dfrac1{(cz+d)^2}\right)$, $\beta=\beta(z)=\Imm\left(\dfrac1{(cz+d)^2}\right)$, $z=\gamma^{-1}(\tilde x+i\tilde y)$.
Pick any observable $\tilde a\in C^\infty_0(T^*\HH)$. Consider operator $\Op_\hbar\tilde a$. Pull-back of such operator by $\gamma$  is $\Op_\hbar a+O_{L^2_{\loc}(\HH)\to L^2_{\loc}(\HH)}(\hbar)$ with 
$$
a(x+iy,\xi_1,\xi_2) = \tilde a\left(\gamma(x+iy),\dfrac{\alpha \xi_1-\beta\xi_2}{\alpha^2+\beta^2},\dfrac{\beta \xi_1+\alpha\xi_2}{\alpha^2+\beta^2}\right), ~~~ ({x},y, {\xi_1}, \xi_2) \in T^*\HH.
$$
This means that if $V(z) = (\Op_{\hbar} \tilde a) v_\hbar|_{\gamma z}$ then also $V(z)=\Op_{\hbar} a (v_\hbar\circ\gamma)|_{z}$, up to minor corrections in $L^2_{\loc}(\HH)$. Thus we have:
\begin{multline}
\label{eq:qform_push_forward}\left\langle\left(\Op_\hbar\tilde a\right) v_\hbar, v_\hbar\right\rangle=
\int_{\HH} \left[\left(\Op_\hbar\tilde a(\tilde z,\tilde\xi)\right) v_\hbar\right](\tilde z)\cdot \bar v_\hbar(\tilde z)\,d\mathcal{A}(\tilde z)=\int_{\HH}V(\gamma^{-1}\tilde z)\cdot \bar v_\hbar(\tilde z)\,d\mathcal A(\tilde z)=\\=\int_\HH V(z)\cdot \bar v_\hbar(\gamma z)\,d\mathcal A(z)=
\int_{\HH}  \left[\Op_\hbar a(z,\xi) (v_\hbar\circ\gamma) \right](z)\cdot \bar v_\hbar(\gamma z)\,d\mathcal{A}(z)+O(\hbar)=\\=
\int_{\HH} \left[\left(\frac{cz+d}{c\bar z+d}\right)^{-\tk/\hbar}\cdot \Op_\hbar a(z,\xi)\cdot \left(\frac{cz+d}{c\bar z+d}\right)^{\tk/\hbar}v_\hbar(z)\right]\cdot \bar v_\hbar(z)\,d\mathcal{A}(z)+O(\hbar).
\end{multline}
We used the fact that  $z\mapsto\bar{z}=\gamma z$ is an isometric change of variable;  the last relation is true because $v_\hbar$ belongs to $\Sp^{\tk/\hbar}(\Gamma)$ and possesses the corresponding twisted automorphy.

In (\ref{eq:qform_push_forward}), we arrived to $\left(\dfrac{cz+d}{c\bar z+d}\right)^{-\tk/\hbar}\cdot \Op_\hbar a(z,\xi)\cdot \left(\dfrac{cz+d}{c\bar z+d}\right)^{\tk/\hbar}$. Here, rational factors  do not almost commute with the central one. We deal with this product using  Yu.~Egorov Theorem. 

To this end, write $\dfrac{cz+d}{c\bar z+d}$ as $\exp(-ip(z,\xi))$, where  phase $p(z,\xi)\colon T^*(\HH)\to\R$ does not depend on $\xi$ and is defined as $p(z,\xi):=-2\arg(cz+d)=2\arctg\left(\dfrac{cx+d}{cy}\right)$. We have $\left(\dfrac{cz+d}{c\bar z+d}\right)^{\tk/\hbar}=\exp\left(\dfrac{\tk\Op_\hbar p}{i\hbar}\right)$. Let $\exp t\Xi_p\colon T^*\HH\to T^*\HH$, $t\in\R$, be Hamiltonian flow defined by $p$. 

We use Yu. Egorov Theorem in the form given in \cite{DimassiSjostrand}, \cite{Zw}. (Notice that this differs from the original result in \cite{Egorov} which is more general but is not implemented in $\hbar$-pseudodifferential operators.) So, by this Theorem we have 
\begin{equation}
\label{eq:p_hamiltonian}
\left(\frac{cz+d}{c\bar z+d}\right)^{-\tk/\hbar}\cdot \Op_\hbar a \cdot \left(\frac{cz+d}{c\bar z+d}\right)^{\tk/\hbar} = \Op_\hbar(a\circ \exp \tk\Xi_p)+O_{L^2_{\loc}(\HH)\to L^2_{\loc}(\HH)}(\hbar).
\end{equation}
Differential equations for $\exp t\Xi_p$ are
$$
\dot x=0, ~ \dot y=0, ~ \dot\xi_1=-\frac{2c^2y}{(cx+d)^2+(cy)^2},~\dot{\xi}_2= \frac{2c(cx+d)}{(cx+d)^2+(cy)^2}.
$$
Thus, 
\begin{multline}
\label{eq:symbol_change}
(a\circ \exp \tk\Xi_p)(x+iy,\xi_1,\xi_2) = a\left(x+iy,\xi_1-\tk\cdot\frac{2c^2y}{(cx+d)^2+(cy)^2},\xi_2+\tk\cdot\frac{2c(cx+d)}{(cx+d)^2+(cy)^2}\right)=\\
=\tilde a\left(\gamma(x+iy),\dfrac{\alpha \xi_1-\beta\xi_2}{\alpha^2+\beta^2}+2\tk c^2y,\dfrac{\beta \xi_1+\alpha\xi_2}{\alpha^2+\beta^2}+
2\tk c(cx+d)\right).
\end{multline}
So, if we put 
$$
\delta(x+iy,\xi_1,\xi_2):=\left(\gamma(x+iy),\dfrac{\alpha \xi_1-\beta\xi_2}{\alpha^2+\beta^2}+2\tk c^2y,\dfrac{\beta \xi_1+\alpha\xi_2}{\alpha^2+\beta^2}+
2\tk c(cx+d)\right),
$$
then, by (\ref{eq:qform_push_forward}), (\ref{eq:p_hamiltonian}) and (\ref{eq:symbol_change}), we have 
$$
\int_{T^*\HH} \tilde a\,d\bar\mu = \int_{T^*\HH} (\tilde a\circ\delta)\,d\bar\mu,
$$
and if $b\in C_0^\infty(T\HH)$ is some function then 
$$
\int_{T\HH} b\,d\underline\mu=\int_{T^*\HH} (b\circ\phi_\tk)\,d\bar\mu=\int_{T^*\HH} (b\circ\phi_\tk\circ\delta)\,d\bar\mu=\int_{T\HH} (b\circ\phi_\tk\circ\delta\circ\phi_\tk^{-1})\,d\underline\mu.
$$
So, in order to prove that $\underline{\mu}=(\phi_{\tk})_{\sharp}{\bar\mu}$ is invariant under action of $\Gamma$, we just have to check that 
$\phi_\tk\circ\delta\circ\phi_{\tk}^{-1}$ is the differential of $\gamma$ given by
$$
(x+iy,v_x,v_y)\mapsto(\gamma(x+iy),\alpha v_x-\beta v_y,\beta v_x+\alpha v_y).
$$
But this can be done by a  straightforward computation. $\blacksquare$

\medskip

\pagebreak

\noindent Now we may also conclude the 

\medskip

\noindent {\bf Proof of Proposition \ref{prop:invariance}}. First assertion is an immediate consequence of Lemma~\ref{lemma:invariance}.

For the second one, take any two locally isometric coverings $\pi_{\rm I},\pi_{\rm II}\colon \HH\to X$. Then there exists hyperbolic isometry $\gamma$ of the form $\gamma z = \dfrac{az+b}{cz+d}$ such that $\pi_{\rm I} = \pi_{\rm II}\circ \gamma$. Function $u_0$ was initially defined on $X$; put $u_{j}:=u_0\circ\pi_j$, $j=\rm I, \rm II$. Then $u_{\rm I} = u_{\rm II}\circ\gamma$.

It is known that for any $f\in C^\infty(\HH)$ and any $\gamma\in\Isom^+(\HH)$ with $\gamma z = \dfrac{az+b}{cz+d}$, and for any $\tau\in\mathbb Z$, we have 
$$
K_\tau\left[f(\gamma z)\left(\dfrac{c\bar z+d}{cz+d}\right)^{\tau}\right]=
\left(\dfrac{c\bar z+d}{cz+d}\right)^{\tau+1}(K_\tau f)|_{\gamma z}
$$
(see \cite{Fay}). Then, by induction,
$$
K_{\tau-1} \dots K_1 K_0 u_{\rm I}(z) = \left(\dfrac{c\bar z+d}{cz+d}\right)^{\tau}\cdot K_{\tau-1} \dots  K_1 K_0 u_{\rm II}|_{\gamma z} 
$$
for any integer $\tau$, and if
$u^{(s)}_{j,\tk} = \Exc^{(s)}_{0\to \tk} u_{j}$, $j=\rm I, \rm II$, then 
\begin{equation}
\label{eq:u_spin}
u^{(s)}_{{\rm II},\tk}(\gamma z) = \left(\dfrac{cz+d}{c\bar z+d}\right)^{[\tk s]} u^{(s)}_{\rm I, \tk}(z).
\end{equation}
Let $\bar{\mu}^\tk_{j}$, $j=\rm I, \rm II$, be semiclassical measures for sequences $\{(u_{j,\tk}^{(s)},1/s)\}$, $s$ ranging some discrete set accumulating to $\infty$. Exactly the same computation as in the proof of Lemma \ref{lemma:invariance} leads us to the  relation
$$
(\phi_\tk)_\sharp \bar{\mu}^\tk_{\rm II} = (D\gamma)_{\sharp}(\phi_\tk)_\sharp \bar{\mu}^\tk_{\rm I}.
$$
Indeed, in that proof we have not used group property and relied only on (\ref{eq:v_spin})  which is the same as (\ref{eq:u_spin}). But this means that the semiclassical measure transferred from cotangent bundle to tangent bundle does not depend on the covering. Thus, the second assertion is also proven.
$\blacksquare$

\section{Ascension during finite adiabatic time}
\label{sec:finite_time}

In this Section we are going to prove Theorem \ref{th:finite_time_shift}. 

First, we cover $X$ by hyperbolic cylinder $\Cyl_l$ with neck length $l>0$ which is the surface $\langle z\mapsto e^lz\rangle\setminus\mathbb C^+$, the Lobachevsky plane folded by the cyclic group of its hyperbolic isometries spanned by transformation $z\mapsto e^lz$, $z\in\mathbb C$. Usual theory of coverings allows to construct a plenty of cylindric covers of~$X$.

Introduce a coordinate system on $\Cyl_l$. First, write $z=x+iy\in\mathbb C^+$ as $z=i\exp(\sigma-i\beta)$, $\sigma\in\R$, $\beta\in\left(-\dfrac{\pi}2,\dfrac{\pi}2\right)$. Then $(\beta,\sigma)$ gives a good and conformal coordinate system on $\Cyl_l$ (now $\sigma$ ranges $\R$ modulo $l$).  Let $\xi$  be coordinate conjugate to $\beta$, and $\eta$ be coordinate conjugate to $\sigma$,  so that $(\beta,\sigma,\xi,\eta)$ is a canonical coordinate system in $T^*\Cyl_l$.  

Operators $K_\tau$, $\tau\in\mathbb Z$, commute with change variable $z\mapsto e^l z$. Thus, all the $K_\tau$ are well-defined on $\Cyl_l$. The same is true for all the $D^\tau$.  We identify function $u_0$ on $X$ with its lift on $\Cyl_l$. Functions $u_\tk$ for $\tk\in\R$ are also defined on cylinder $\Cyl_l$ since operators $K_\tau$ are. 

Magnetic Hamiltonian  on $\Cyl_l$ takes the form 
$$
H_\tk = \frac{(\xi-\tk)^2\cos^2\beta +(\eta\cos\beta-\tk\sin\beta)^2}2.
$$
Pick positive $\emm$ small enough such that 
\begin{equation}
\label{eq:eta_max_def}
\emm < 1/2, ~~ \tk\emm < \sqrt{1-\emm^2}.
\end{equation}
Put 
\begin{equation}
\label{eq:Omega}
\Omega^0 := \{(\beta,\sigma,\xi,\eta)\in T^*\Cyl_l\colon |\eta|<\emm, \, \xi >0, \, H_0(\beta,\sigma,\xi,\eta) = 1/2\}
\end{equation}
and
\begin{equation}
\Omega^\tk := \{(\beta,\sigma,\xi,\eta)\in T^*\Cyl_l\colon |\eta|<\emm, \, \xi >0, \, H_\tk(\beta,\sigma,\xi,\eta) = (\tk^2+1)/2\}.
\end{equation}
First, we prove the following 

\begin{predl}
\label{prop:cyl_mapping_exists}
There exist a smooth mapping $G\colon \Omega^0\to\Omega^\tk$ \emph{(}in fact, a diffeomorphism\emph{)} and also a smooth function $A\colon \Omega^\tk\to(0,+\infty)$ such that $\bar\mu^\tk\cdot\mathds 1_{\Omega^\tk} = A\cdot G_\sharp (\bar\mu^0\cdot\mathds 1_{\Omega^0})$.
Mapping $G$ and function $A$ are given in a form which does not depend on the initial sequence~$\{u_0^{(s)}\}_{s\in J}$.
\end{predl}

\noindent But these $G$ and $A$ are rather implicit, so the second step is to test the transformations
\begin{equation*}
\overline T\colon \Meas(\Omega^0)\to\Meas(\Omega^\tk), ~~ 
\overline T\mu:= A\cdot G_\sharp\mu, ~~ \mu \in \Meas(\Omega^0),
\end{equation*}
$$
\underline T := (\phi_{\tk})_\sharp\overline T(\phi_{0}^{-1})_\sharp\colon \Meas\left(\phi_0(\Omega^0)\right)\to\Meas\left(\phi_\tk(\Omega^\tk)\right)
$$
by inserting to $\underline T$ a measure concentrated on a single geodesic line. This will allow us to replace calculation of integrals and implicit functions by evaluation of couple of asymptotics. This is done in Subsection~\ref{subsec:testing}.

Before we pass to proofs, let us give some  empiric observations.

\subsection{Numerical experiment: intuition of traveling waves}

\label{subsec:experiment}
This Subsection is mostly informal; we want to  clarify what is happening when we apply $K_\tau/\sqrt{s^2+\tau(\tau+1)}$ to an eigenfunction $u$ with $D^\tau u = s^2 u$.  

Let's separate variables in operator $D^{\tau}$. A good way is to search for eigenfunctions of the form $u=\exp(iax)w(y)$, where $(x,y)$ are standard coordinates in $\mathbb C^+$ and $a\in\R$ is a parameter. We put $s_1 := \sqrt{s^2-1/4}$ so that $s^2=s_1^2+1/4$. Equation 
$$
D^\tau \left(\exp(iax)w(y)\right)=\left(s_1^2+1/4\right)\exp(iax)w(y)
$$
then becomes
\begin{equation}
\label{eq:Whittaker_ODE}
w''(y)+\left(-a^2+\frac{2\tau a}{y}+\frac{s_1^2+1/4}{y^2}\right)w(y)=0.
\end{equation}
This ODE has two linearly independent solutions, one of them is $W_{\tau,is_1}(2ay)$ for $a>0$ (and $W_{-\tau,is_1}(2|a|y)$ for $a<0$), the Whittaker $W$-function. See \cite{Buchholz} about Whittaker functions (our function is $W_{\tau,is_1/2}(2|a|y)$ in the notation of the Buchholz's treatise \cite{Buchholz}). 

We just check the two solutions mentioned above. It seems that any function $u\colon \HH\to \mathbb C$ with $D^\tau u=(s_1^2+1/4) u$ \emph{having tempered growth} can be decomposed into combination of these $W$-functions; the other solution of (\ref{eq:Whittaker_ODE}), Whittaker $M$-function, has exponential growth at $y\to+\infty$ and seems to be unable to contribute to, say, a bounded eigenfunction $u$. Thus, we are going to study ascension evolution of a $W\mbox{-function}$. Since we are in an heuristic considerations, we restrict ourselves to the case $a>0$.

To calculate the derivative, we need a \emph{contiguous relation} on Whittaker functions.  In \cite[p. 81]{Buchholz} we find:
$$
\frac{dW_{\tau,is_1}(y)}{dy}=\left(\frac12-\frac\tau y\right) W_{\tau,is_1}(y)-\frac{W_{\tau+1,is_1}(y)}y.
$$
This yields
$$
\frac{K_\tau}{\sqrt{s_1^2+(\tau+1/2)^2}}(e^{iax} W_{\tau,is_1}(2ay)) = \frac{e^{iax}W_{\tau+1,is_1}(2ay)}{\sqrt{s_1^2+(\tau+1/2)^2}}.
$$

Notice that if $a$ varies, $e^{iax} W_{\tau,is}(2ay)$ stays a rescale of one fixed function; this is because $z\mapsto az$ is an isometry of Lobachevsky plane $\HH\simeq \mathbb C^+$. Thus, we may take arbitrary $a$ to build graphics.

\begin{figure}
	\centering
	{\includegraphics[scale=0.751]
		{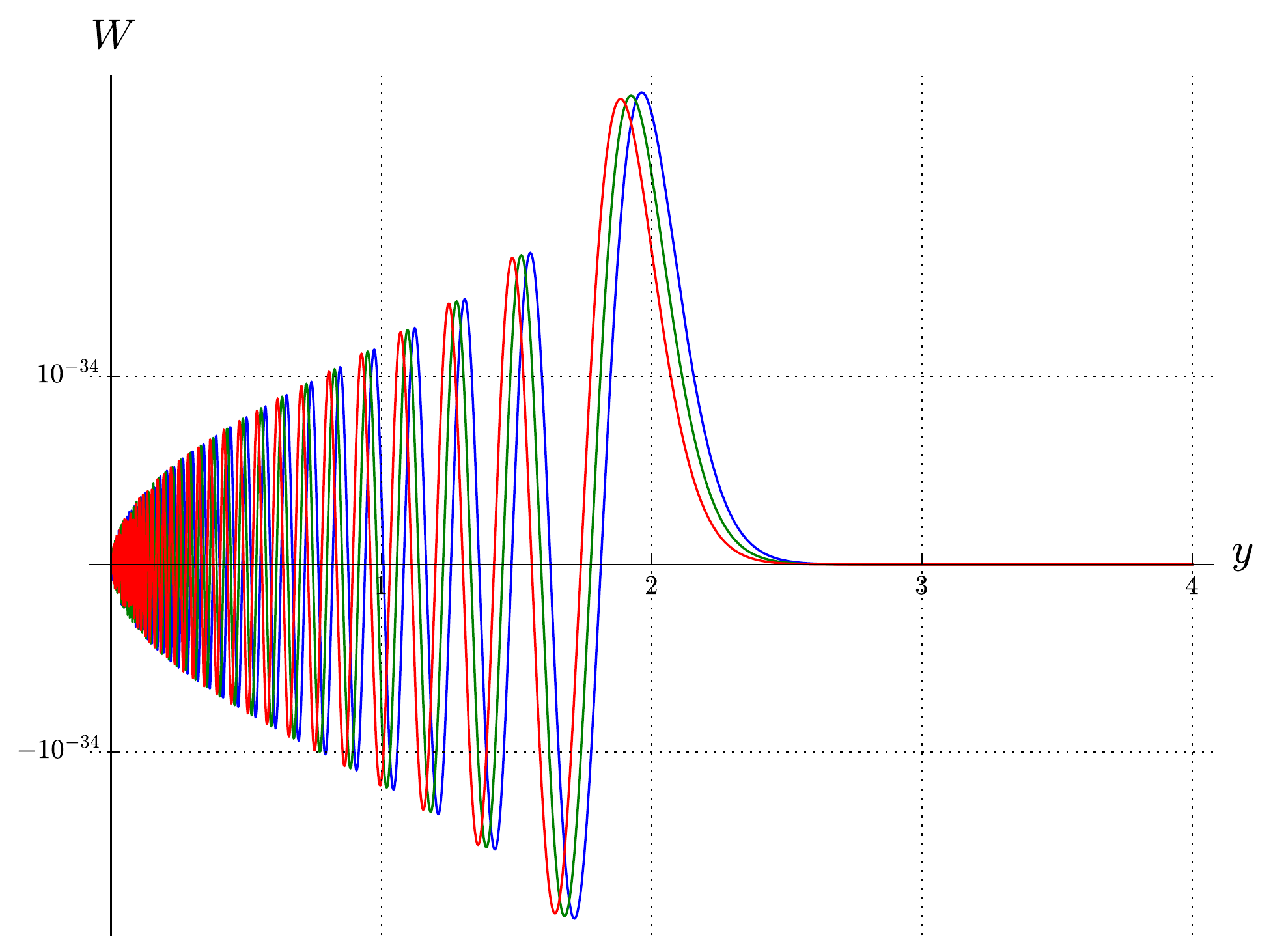}}
	\caption{Ascension evolution of a Whittaker harmonic}
	\label{fig:whittaker_waves}
\end{figure}

We successfully use \cite{Sage} and \cite{mpmath} to calculate hypergeometric functions with large parameters. At Figure \ref{fig:whittaker_waves}, the left-shifted (red) 
wave is $W_{0,\,i\cdot50}(50y)$, the middle one (green) is $W_{1,\,i\cdot50}(50y)/\sqrt{50^2+1/4}$, the one shifted to the right (blue) is $W_{2,\,i\cdot50}(50y)/\sqrt{(50^2+1/4)\cdot(50^2+9/4)}$. We see that the wave $W_{0,\,i\cdot50}(\cdot)$ runs to the right under ascension evolution. Let us evaluate the speed of this running. For that, we find the abscissae and ordinata of the rightest and highest  peaks of these waves: the abscissae are $1.884...,1.922...,1.962...$, the ordinata are $2.488...\cdot 10^{-34},2.499...\cdot 10^{-34},2.510...\cdot 10^{-34}$. We see that abscissae of peaks moved at about $0.04$ which is comparable to $\hbar\approx1/s$. The same will be true for ordinata if we normalize amplitudes of our waves multiplying them by $10^{34}$.

The same behavior of ascension evolution is exhibited for all the other values of $\tau$ and~$s$. To summarize, \emph{$\dfrac{K_\tau(e^{iax}W_{0,is_1}(2ay))}{\sqrt{s_1^2+(\tau+1/2)^2}}$ is the wave $e^{iax} W_{0,is_1}(2ay)$ shifted by a distance not exceeding $\const\cdot\hbar$, $\hbar=1/s$}. Then, \emph{as $\tk\approx\hbar\tau$ varies \emph{(}almost\emph{)} smoothly, the wave $\Maas_{0\to\tk}^{(s)} (e^{iax}W_{0,is_1}(2ay))$  runs with a finite speed and also almost smoothly}.

\medskip

We proceed by an heuristic observation concerning propagation of waves. Suppose that we are given by some group of linear transformations $E^\tk$ (parametrized by $\tk\in\R$) acting on functions on an Euclidean space  $\R^d$, $d\in\mathbb N$, or somewhere else; suppose, further, that there is some set of harmonics $\{w_\eta\}_{\eta\in I}$ in $\R^d$ (the notion of a harmonic is understood in a very wide sense) and that, for any $\eta$, the function $E^\tk w_\eta$ depends on $\tk$ as a wave running with bounded speed (where speed is calculated with respect to parameter~$\tk$). Then $E^\tk$ moves semiclassical measures corresponding to linear combinations of functions $w_\eta$ in a controlled way. In particular, $E^\tk$ preserves microlocal singularities formed by such combinations. 

For quantum Hamiltonian evolutions, this observation is just Yu. Egorov Theorem. An arbitrary evolution in order to be governed by this informal statement should move wavevectors smoothly but, unlike Hamiltonian flow, does not have to preserve phase space volume. Thus, formalization of our heuristics should be a version of Yu.~Egorov Theorem accomplished by the Jacobian of the classical flow with respect to phase space volume.

This heuristic observation leads us to the idea that \emph{we have a good chance to get a control over measures $\bar{\mu}^\tk$ in terms of measure $\bar{\mu}^0$}.

See also \cite{Taylor} for singularities reconstruction in the case of compact manifold and waves running with unbounded speeds. 

\medskip

Now, let us take a numerical experiment for waves on cylinder. To this end, we separate variables $(\beta,\sigma)$. We have
\begin{gather*}
y\frac{\partial}{\partial x} = \sin\beta\cos\beta\frac{\partial}{\partial\sigma}+\cos^2\beta\frac{\partial}{\partial\beta},\\
y\frac{\partial}{\partial y} = \cos^2\beta\frac{\partial}{\partial\sigma}-\sin\beta\cos\beta\frac{\partial}{\partial\beta}.
\end{gather*}
These two vectorfields are invariant with respect to $z\mapsto e^lz$ and thus can be correctly defined on $\Cyl_l$. (In particular, this implies that Maa\ss{} derivatives $K_\tau$ can also be defined on the cylinder.) Let us search for eigenwaves of the form $\exp(im\sigma)\cdot w(\beta)\colon \Cyl_l\to\mathbb C$, $m\in2\pi\mathbb Z/l$. We have  
\begin{multline}
\label{eq:cylinder_raising}
K_\tau\left(\exp(im\sigma)\cdot w(\beta)\right)=\left(2iy\frac{\partial}{\partial z}+\tau\right)\left(\exp(im\sigma)\cdot w(\beta)\right)=\\=
\exp(im\sigma)\cdot \left(\tau w+i\cos\beta\cdot  e^{i\beta}mw+i\cos\beta\cdot e^{i\beta} w'\right),
\end{multline}
\vspace{-1.cm}
\begin{multline}
\label{eq:cylinder_eigen}
\vspace{5cm}D^\tau\left(\exp(im\sigma)\cdot w(\beta)\right)=\left(-\Delta_\HH+2i\tau y\dfrac{\partial}{\partial x}
\right)\left(\exp(im\sigma)\cdot w(\beta)\right) =\\=
\exp(im\sigma)\cdot\left(-\cos^2\beta\cdot w''+2i\tau\cos^2\beta\cdot w'+\left(m^2\cos^2\beta-2\tau m \sin\beta\cos\beta\right)w\right).
\end{multline}

We can find solutions of $\exp(-im\sigma)D^\tau\left(\exp(im\sigma)\cdot w(\beta)\right)=(s_1^2+1/4)w(\beta)$ explicitly in terms of hypergeometric functions:
\begin{gather*}
w_{{\rm I},\tau,m}(\beta)=\cos^{-\tau}\beta\cdot e^{\beta(-m+i\tau)} {}_2F_1\left(\frac12+is_1+\tau, \frac12-is_1+\tau;-im+\tau+1;\frac12-i\frac{\tg\beta}2\right),\\
w_{{\rm II},\tau,m}(\beta)=w_{{\rm I},-\tau,m}(-\beta).
\end{gather*}
One has
$$
\frac{K_\tau\left(\exp(im\sigma)w_{{\rm I},\tau,m}(\beta)\right)}{\sqrt{s_1^2+\left(\tau+\frac12\right)^2}} =
\exp(im\sigma) w_{{\rm I},\tau+1,m}(\beta)\cdot\frac{\sqrt{s_1^2+\left(\frac12+\tau\right)^2}}{2(1+\tau-im)},
$$
and the same for the second solution.
These expressions allow us to perform numerical experiments as well as above. We again observe the same behavior: these waves travel with bounded speed with respect to adiabatic time parameter $\tk=\hbar\tau$, $\hbar=1/s$.

We make one more numerical observation useful in the below. Unlike Whittaker functions from the beginning of this Subsection, the above waves $w_{{\rm I},\tau,m}(\beta), w_{{\rm II},\tau,m}(\beta)$ are \emph{monochromatic} in the sense that their absolute values do not oscillate (this is an experimental fact); this means that each of the wave has only one microlocal frequency. Hence, \emph{ascension evolution preserves monochromaticity of cylindrical harmonics}. This observation will allow us to distinguish between the two solutions of eigenwave equation in Lemma \ref{lemma:monochrome} in the following Subsection.

\subsection{WKB ansatz for cylindrical harmonics}
\label{subsec:WKB}
Now, let us focus on rigorous analytical study of cylindrical harmonics. Recall that we deal with $u_{\tk} = \Exc_{0\to \tk}^{(s)} u_0$, this function satisfies $D^\tau u_{\tk}= s^2 u_{\tk}$ with $\tau=[\tk s]$. So, in this Subsection we apply the very standard WKB techniques to study the equation
\begin{equation}
\label{eq:cyl_harm} \exp(-im\sigma)D^\tau\left(\exp(im\sigma)\cdot w(\beta)\right)=s^2w(\beta)
\end{equation}
on cylindric harmonic. We mostly calculate two higher order  asymptotic terms. Put  $\tm=m/s$ and $\tk_1=\tau/s$. The latter $\tk_1$ is a temporary denotation, it will be used only in this Subsection and also in the next one. Note that $\tm$ is the same as $\eta$ but ranges a discrete set.

Using (\ref{eq:cylinder_eigen}), we rewrite equation (\ref{eq:cyl_harm}) in $Q$-form:
$$
\left(w(\beta)\cdot e^{-i\tau\beta}\right)''+s^2\cdot Q_{\tk_1,\tm}(\beta)\cdot w(\beta)e^{-i\tau\beta}=0,
$$
where $Q=Q_{\tk_1,\tm}(\beta)=2\tk_1\tm\tg\beta-\tm^2+\dfrac{1}{\cos^2\beta}+\tk_1^2$. Therefore, (\ref{eq:cyl_harm}) has two WKB solutions satisfying the asymptotics
\begin{equation}
\label{eq:WKB_main}
w^{{\rm I},s}_{\tk_1,\tm}(\beta), \, w^{{\rm II},s}_{\tk_1,\tm}(\beta)=\exp\left(i\tau\beta\pm is\cdot\int\limits_0^\beta \sqrt{Q_{\tk_1,\tm}(\beta_1)}\,d\beta_1 -\frac14\int\limits_0^\beta\frac{Q'_{\tk_1,\tm}(\beta_1)}{Q_{\tk_1,\tm}(\beta_1)}\,d\beta_1+O\left(s^{-1}\right)\right)
\end{equation}
(see \cite{Fedoryuk}).  Constant implied in remainder $O(s^{-1})$ is uniform over compact sets of parameters $(\beta,\tk_1,\tm)$ ranging away from \emph{turning points}, that is, from zeroes of $Q_{\tk_1,\tm}(\beta)$. Automatically, this holds when $\tk_1$ ranges any compact interval in $\R$, $|\tm|\le1/2$ and $\beta$ ranges any compact subinterval in $\left(-\dfrac{\pi}2,\dfrac{\pi}2\right)$ --- since $Q_{\tk_1,\tm}(\beta)$ does not vanish therein.  

We will mostly study the first WKB solution, $w^{{\rm I},s}_{\tk_1,\tm}(\beta)$, the one with $+$ sign before the imaginary exponent in (\ref{eq:WKB_main}). The second solution will be eliminated from the expansion of a general eigenfunction by a corresponding frequency cut-off in Subsection~\ref{subsubsec:cutoff}.

Now, we formalize the empiric observation from the end of Subsection \ref{subsec:experiment}.

\begin{lemma}
	\label{lemma:monochrome}
	For WKB solutions defined in \emph{(\ref{eq:WKB_main})} above and $\tau=\tk_1 s$, we have
	\begin{multline}
	\label{eq:maas_diff_ansatz}
	e^{-im\sigma}\frac{K_\tau}{\sqrt{s^2+\tau(\tau+1)}}\left(e^{im\sigma} w^{{\rm I},s}_{\tk_1,\tm}(\beta)\right)=\\=\left(c_1(\tk_1,\tm,s)+O\left(s^{-2}\right)\right)\cdot w^{{\rm I},s}_{\tk_1+1/s,\tm}(\beta)+O(s^{-2})\cdot w^{{\rm II},s}_{\tk_1+1/s,\tm}(\beta),
	\end{multline}
	where
	\begin{multline*}
	c_1 = c_1(\tk_1,\tm,s) = \frac{\sqrt{\tk_1^2-\tm^2+1}-i\tm}{\sqrt{1+\tk_1^2}}+\\+\frac{1}{s}\cdot\left(\frac{\left(-\sqrt{\tk_1^2-\tm^2+1}+i\tm\right)\tk_1}
	{2\left(1+\tk_1^2\right)}+\frac{i\tm\tk_1}{2\left(\tk_1^2-\tm^2+1\right)}\right)\cdot\frac{1}{\sqrt{\tk_1^2+1}}.
	\end{multline*}
	The error terms $O(s^{-2})$ are both complex scalars uniform when $\tk_1$ ranges any compact interval in $\R$, $|\tm|\le 1/2$ and $\beta$ ranges any compact subset in $\left(-\dfrac{\pi}2,\dfrac{\pi}2\right)$.
	
	In other words, ascension evolution \emph{preserves monochromaticity} of the $w^{\rm I}$-wave up to minor corrections.
\end{lemma}

\noindent {\bf Remark.} The same statement holds for $w^{\rm II}$-waves. Namely, if we put
\begin{multline*}
c_2^{\rm II}(\tk_1,\tm,s) = \frac{-\sqrt{\tk_1^2-\tm^2+1}-i\tm}{\sqrt{1+\tk_1^2}}+\\+\frac{1}{s}\cdot\left(\frac{\left(\sqrt{\tk_1^2-\tm^2+1}+i\tm\right)\tk_1}
{2\left(1+\tk_1^2\right)}+\frac{i\tm\tk_1}{2\left(\tk_1^2-\tm^2+1\right)}\right)\cdot\frac{1}{\sqrt{\tk_1^2+1}}
\end{multline*}
then we have 
\begin{equation*}
e^{-im\sigma}\frac{K_\tau}{\sqrt{s^2+\tau(\tau+1)}}\left(e^{im\sigma} w^{{\rm II},s}_{\tk_1,\tm}\right)=O(s^{-2})\cdot w^{{\rm I},s}_{\tk_1+1/s,\tm}+\left(c_2^{\rm II}(\tk_1,\tm,s)+O(s^{-2})\right)\cdot w^{{\rm II},s}_{\tk_1+1/s,\tm}
\end{equation*}
(we just changed signs at $\sqrt{\tk_1^2-\tm^2+1}$).

\medskip

\noindent {\bf Proof of Lemma \ref{lemma:monochrome}.} Since 
$$
D^{\tau}\left(\exp(im\sigma)w^{{\rm I},s}_{\tk_1,\tm}(\beta)\right)=s^2\cdot \exp(im\sigma)w^{{\rm I},s}_{\tk_1,\tm}(\beta),
$$ 
by Proposition \ref{prop:Maas_properties} we have 
$$
D^{\tau+1}\frac{K_\tau}{\sqrt{s^2+\tau(\tau+1)}}\left(\exp(im\sigma)w^{{\rm I},s}_{\tk_1,\tm}(\beta)\right)=s^2\cdot\frac{K_\tau}{\sqrt{s^2+\tau(\tau+1)}}\left( \exp(im\sigma)w^{{\rm I},s}_{\tk_1,\tm}(\beta)\right).
$$ 
This leads to 
\begin{equation*}
e^{-im\sigma}\frac{K_\tau}{\sqrt{s^2+\tau(\tau+1)}}\left(e^{im\sigma} w^{{\rm I},s}_{\tk_1,\tm}\right)=c_1\cdot w^{{\rm I},s}_{\tk_1+1/s,\tm}+c_2\cdot w^{{\rm II},s}_{\tk_1+1/s,\tm},
\end{equation*} with some scalars $c_1, c_2\in\mathbb C$. Having (\ref{eq:cylinder_raising}) in mind, put 
\begin{gather*}
v(\beta):=s\tk_1 w^{{\rm I},s}_{\tk_1,\tm}+is\cos\beta\cdot  e^{i\beta}\tm w^{{\rm I},s}_{\tk_1,\tm}+i\cos\beta\cdot e^{i\beta}\cdot \frac{dw^{{\rm I},s }_{\tk_1,\tm}}{d\beta},\\
w(\beta) := \frac{v(\beta)}{\sqrt{s^2+s\tk_1\cdot(s\tk_1+1)}}.
\end{gather*}
To find $c_1,c_2$,  we consider asymptotics in $\beta=0$ of all the $w$-functions and write the $2\times 2$ system on $c_1, c_2$ as
\begin{equation}
\hspace{-1cm}
\left\{\begin{aligned}
c_1\cdot w^{{\rm I},s}_{\tk_1+1/s,\tm}(0)+c_2\cdot w^{{\rm II},s}_{\tk_1+1/s,\tm}(0)=& w(0)\\
c_1\cdot\frac1s \cdot\frac{dw^{{\rm I},s}_{\tk_1+1/s,\tm}}{d\beta}(0)+c_2\cdot\frac1s\cdot\frac{dw^{{\rm II},s}_{\tk_1+1/s,\tm}}{d\beta}(0)=& \frac1s\cdot\left.\frac{dw}{d\beta}\right|_{\beta=0}
\end{aligned}
\right.
\label{eq:diff_sys2x2}
\end{equation}

From the construction of WKB solutions in \cite{Fedoryuk}, it can be seen that we may assume \emph{precise} equalities
\begin{gather*}
\frac{dw^{{\rm I},s}_{\tk_1,\tm}}{d\beta}(0) = is\tk_1+is\sqrt{Q_{\tk_1,\tm}(0)}-\frac{Q_{\tk_1,\tm}'(0)}{4Q_{\tk_1,\tm}(0)},\\
\frac{dw^{{\rm II},s}_{\tk_1,\tm}}{d\beta}(0) = is\tk_1-is\sqrt{Q_{\tk_1,\tm}(0)}-\frac{Q_{\tk_1,\tm}'(0)}{4Q_{\tk_1,\tm}(0)},\\
w^{{\rm I},s}_{\tk_1,\tm}(0)=w^{{\rm II},s}_{\tk_1,\tm}(0) = 1.
\end{gather*}
The similar holds for $\tk_1$ replaced by $\tk_1+1/s$. We see that, for $s$ large, coefficients of the system~(\ref{eq:diff_sys2x2}) are of $O(1)$ order and that its determinant  is separated from zero; in other words, this system is well-posed. Hence, to prove our Lemma, it is enough to check that $c_1,c_2$ from the statement of satisfy system~(\ref{eq:diff_sys2x2}) up to $O(1/s^2)$ errors. We proceed by opening the brackets.

We find $Q_{\tk_1,\tm}(0)=\tk_1^2-\tm^2+1$, $Q_{\tk_1,\tm}'(0)=2\tk_1\tm$. From the equation (\ref{eq:cyl_harm}) on $w^{{\rm I},s}_{\tk_1,\tm}$ and from (\ref{eq:cylinder_eigen}), we have
\begin{multline*}
\left.\frac{d^2w^{{\rm I},s}_{\tk_1,\tm}}{d\beta^2}\right|_{\beta=0}=\left.2is\tk_1 \cdot \frac{dw^{{\rm I},s}_{\tk_1,\tm}}{d\beta}\right|_{\beta=0}+\left.s^2\cdot\left(\tm^2-2\tk_1\tm\tg\beta-\cos^{-2}\beta\right)w^{{\rm I},s}_{\tk_1,\tm}\right|_{\beta=0} = \\=
-s^2\cdot\left(\tk_1^2+2\tk_1\sqrt{Q_{\tk_1,\tm}(0)}+Q_{\tk_1,\tm}(0)\right)-s\cdot\frac{i\tk_1 Q_{\tk_1,\tm}'(0)}{2Q_{\tk_1,\tm}(0)}.
\end{multline*}
Substituting, we find
$$v(0)=s\cdot\left(i\tm-\sqrt{Q_{\tk_1,\tm}(0)}\right)-\frac{i Q_{\tk_1,\tm}'(0)}{4Q_{\tk_1,\tm}(0)},
$$
\begin{multline*}
v'(0)=-s^2\cdot\left(i\tk_1\sqrt{Q_{\tk_1,\tm}(0)}+i Q_{\tk_1,\tm}(0)+\tk_1\tm+\tm\sqrt{Q_{\tk_1,\tm}(0)}\right)+\\+s\cdot\left(\frac{(\tk_1-i\tm)
 Q_{\tk_1,\tm}'(0)}{4Q_{\tk_1,\tm}(0)}-\tm-i\tk_1-i\sqrt{Q_{\tk_1,\tm}(0)}\right).
\end{multline*}
Denominator in $w$ expands as
$$
b:=\frac{1}{\sqrt{s^2+s\tk_1\cdot(s\tk_1+1)}} = \frac{1}{s\sqrt{1+\tk_1^2}}-\frac{\tk_1}{2s^2(1+\tk_1^2)^{3/2}}+O(1/s^3).
$$
Hence, if we put 
\begin{multline}
\label{eq:c1_def}
c_1 := b\cdot v(0)=\\=\left(\frac{1}{s\sqrt{1+\tk_1^2}}-\frac{\tk_1}{2s^2(1+\tk_1^2)^{3/2}}\right)\cdot \left(s\cdot\left(i\tm-\sqrt{Q_{\tk_1,\tm}(0)}\right)-\frac{i Q_{\tk_1,\tm}'(0)}{4Q_{\tk_1,\tm}(0)}\right)+O(1/s^2)
\end{multline}
and $c_2 := 0$, then these $c_1, c_2$ will enjoy the first equation in (\ref{eq:diff_sys2x2}) up to $O(1/s^2)$. We easily see that $c_1$ from (\ref{eq:c1_def}) is the same as $c_1$ from the statement of Lemma.

By Taylor formula, 
$$
\frac{dw^{{\rm I},s}_{\tk_1+1/s,\tm}}{d\beta}(0) = s\cdot\left(i\tk_1+i\sqrt{Q_{\tk_1,\tm}(0)}\right)+\left(i+\frac{i\tk_1}{\sqrt{Q_{\tk_1,\tm}(0)}}-\frac{Q_{\tk_1,\tm}'(0)}{4Q_{\tk_1,\tm}(0)}\right)+O(1/s).
$$
To verify the second equation of (\ref{eq:diff_sys2x2}), we have to check that 
$$
b\cdot v(0)\cdot\frac{dw^{{\rm I},s}_{\tk_1+1/s,\tm}}{d\beta}(0) = b\cdot v'(0)+O(1/s)
$$
or that 
$$
v(0)\cdot\frac{dw^{{\rm I},s}_{\tk_1+1/s,\tm}}{d\beta}(0) = v'(0)+O(1).
$$
But this is done by a straightforward opening the brackets. $\blacksquare$

\medskip

\noindent We derive a corollary from our ansatz. Put 
\begin{equation}
\label{eq:omega_def}
\omega_{\tk,m,s}(\beta) := w^{{\rm I},s}_{[\tk s]/s,m/s}(\beta)\cdot \prod\limits_{\tau=0}^{[\tk s]-1} c_1(\tau/s,\tm,s).
\end{equation}

\begin{sled}
\label{sled:exc_coeff}	One has
	\begin{equation*}
	\Maas_{0\to \tk}^{(s)} (e^{im\sigma}\cdot w^{{\rm I},s}_{0,\tm}(\beta)) = e^{im\sigma}\cdot \omega_{\tk,m,s}(\beta)+O(1/s),
	\end{equation*}
	the remainder estimate is uniform when $\tk\in\R$, $|\tm| = |m/s| \le 1/2$ and $\beta$ ranges any compact subinterval in $\left(-\dfrac{\pi}2,  \dfrac{\pi}2\right)$.
\end{sled}

\noindent Of course, a similar statement holds for ascension  evolution of $w^{\rm II}$-waves. The proof is easily obtained by rewriting the action of normed raising operator  via multiplication over matrices $2\times 2$.

\subsection{One-dimensional phase transport}
\label{subsec:phase_transport}

So, we reduced the action of ascension  evolution on basic $w^{\rm I}$-waves to relation from Corollary \ref{sled:exc_coeff}. Now let us clarify the notion of "wave traveling with finite speed" mentioned in Introduction and in Subsection \ref{subsec:experiment}. In one-dimensional space (spanned on {$\beta$-direction}) this can be done by setting the correspondence between points having \emph{the same phase} at different moments of adiabatic time $\tk$. This is formalized by the following 

\begin{lemma}
	\label{lemma:Phi_diffeo}
	There exist a smooth mapping $\Phi_{\tk,\tm}(\beta) \colon \left(-\dfrac{\pi}2,\dfrac{\pi}2\right) \to  \left(-\dfrac{\pi}2,\dfrac{\pi}2\right)$ depending smoothly on $\tk\in\R$ and $\tm\in \left[-\dfrac{1}{2},\dfrac 1{2}\right]$ and also smooth \emph{real-valued} scalar functions 
	$$
	f_1(\tk,\beta,\tm), \, f_2(\tk,\beta,\tm), \, f_3(\tk,\beta,\tm)\colon \R\times\left(-\dfrac{\pi}2,  \dfrac{\pi}2\right)\times \left[-\dfrac{1}{2},\dfrac 1{2}\right] \to\R
	$$ 
	such that:
	\begin{enumerate}
		\item for any $\tk\in\R$ and $\tm \in\left[-\dfrac{1}{2},\dfrac 1{2}\right]$, the mapping $\beta \mapsto \Phi_{\tk,\tm}(\beta)$ is a smooth increasing diffeomorphism  of interval $\left(-\dfrac{\pi}2,\dfrac{\pi}2\right)$ onto itself;
		
		\item if we put $f_0^{(s)}(\tk,\beta,\tm):=f_1(\tk,\beta,\tm)+\{\tk s\}\cdot f_2(\tk,\beta,\tm)$ \emph{(}$\{\cdot\}$ being fractional part\emph{)} then we have 
		\begin{equation*}
		  \omega_{\tk,m,s}\left(\Phi_{\tk,\tm}(\beta)\right)=
		  w^{{\rm I},s}_{0,\tm}(\beta)\cdot \exp\left({i f_0^{(s)}(\tk,\beta,\tm)+f_3(\tk,\beta,\tm)+O(1/s)}\right), 
		\end{equation*}
		where the constant in $O(1/s)$ remainder term is uniform when $\tm\in\left[-\dfrac{1}{2},\dfrac 1{2}\right]$ and  $\beta$ ranges any compact subinterval in $\left(-\dfrac{\pi}2,\dfrac{\pi}2\right)$. Recall that~$\omega$ is defined in~\emph{(\ref{eq:omega_def})}.
	\end{enumerate} 
\end{lemma}

\noindent This Lemma says that, under ascension evolution, phase of wave is \emph{transported} smoothly with respect to $\beta$ and $\tm$ (up to minor errors; recall here that the main term of phase of $w^{{\rm I},s}_{\tk_1,\tm}$ is of order~$s$, see ansatz~(\ref{eq:WKB_main})).  Remark also that we explicated the dependence of $f^{(s)}_0$ on $s$ just to be calm. In what follows we will differentiate $f^{(s)}_0$ over $\tm$ and we need just the corresponding smoothness.

\medskip

\noindent {\bf Proof of Lemma \ref{lemma:Phi_diffeo}.} For $\tau=0,1,\dots,[Bs]-1$, we set $\tk_1 := \tau/s$ and apply Lemma~\ref{lemma:monochrome}. Since $|c_1(\tk_1,\tm,s)|=1+O(1/s)$ for $c_1$ from that Lemma, we may write this coefficient as
$$
c_1(\tk_1,\tm,s) = \exp\left(i b_1(\tk_1,\tm)+\frac{1}{s}\left(i b_2(\tk_1,\tm)+b_3(\tk_1,\tm)\right)+O(1/s^2)\right);
$$
here and in the rest of the proof, $b_1, b_2, \dots$ are smooth \emph{real-valued} functions of their arguments.

By Euler--Maclaurin formula, 
\begin{multline*}
\sum_{\tau=0}^{[\tk s]-1} b_1(\tau/s,\tm)=
\int\limits_0^{[\tk s]} b_1(\tau/s,\tm)\,d\tau+\frac{b_1(0,\tm)-b_1\left(\dfrac{[\tk s]}{s},\tm\right)}{2}+O(1/s)=\\=
s\cdot \int\limits_0^\tk b_1(\tk_2,\tm)\,d\tk_2-\{\tk s\}\cdot b_1(\tk,\tm)+\frac{b_1(0,\tm)-b_1\left(\tk,\tm\right)}{2}+O(1/s).
\end{multline*}
Also, 
\begin{equation*}
\sum_{\tau=0}^{[\tk s]-1}\frac{i b_2(\tau/s,\tm)+b_3(\tau/s,\tm)}{s} 
= \int\limits_{0}^\tk\left(i b_2(\tk_2,\tm)+b_3(\tk_2,\tm)\right)\,d\tk_2+O(1/s).
\end{equation*}
We therefore have
\begin{equation*}
\prod\limits_{\tau=0}^{[\tk s]-1} c_1(\tau/s,\tm,s) = \exp\left(is\cdot b_4(\tk,\tm)+i\cdot (b_5(\tk,\tm)+\{\tk s\}\cdot b_6(\tk,\tm))+b_7(\tk,\tm)+O(1/s)\right);
\end{equation*}
note that $b_4(\tk,\tm)=\int\limits_0^\tk b_1(\tk_2,\tm)\,d\tk_2$, $b_1(\tk_2,\tm) = -\arctg\left(\tm/\sqrt{\tk_2^2-\tm^2+1}\right)$, we will use this in the proof of Proposition~\ref{prop:bounded_dist} below.

Now, for any $\tk_1\in\R$, define phase function as
$$
P_{\tk_1,\tm}(\beta) := \tk_1\beta+\int\limits_0^\beta\sqrt{Q_{\tk_1,\tm}(\beta_1)}\,d\beta_1,
$$
it does not depend on $s$ at all; this is the main term of phase of $w_{\tk_1,\tm}^{{\rm I},s}$, see (\ref{eq:WKB_main}). So we may take $\arg w_{\tk_1,\tm}^{{\rm I},s}(\beta) = s\cdot P_{\tk_1,\tm}(\beta) + O(s^{-1})$. 

Note also that for any $\tk_1$ the phase mapping $\beta \mapsto P_{\tk_1,\tm}(\beta)$ is a smooth diffeomorphism of $(-\pi/2,\pi/2)$ onto $(-\infty,+\infty)$; it also depends on $\tm$ smoothly (if $|\tm|\le 1/2$). Thus, if we define $\Phi_{\tk,\tm}(\beta)$, $\beta\in(-\pi/2,\pi/2)$, by 
\begin{equation}
\label{eq:Phi_define}
P_{\tk,\tm}\left(\Phi_{\tk,\tm}(\beta)\right) = P_{0,\tm}(\beta)-b_4(\tk,\tm),
\end{equation}
then this will equalize the higher-order term of phase up to $O(1)$:
$$
\omega_{\tk,m,s}(\Phi_{\tk,\tm}(\beta)) = w_{0,\tm}^{{\rm I},s}(\beta)\cdot \exp(O(1)).
$$
(see (\ref{eq:omega_def})).

More precisely, we have 

\begin{multline*}
\arg\omega_{\tk,m,s}(\Phi_{\tk,\tm}(\beta))=\arg w^{{\rm I},s}_{[\tk s]/s,\tm}(\Phi_{\tk,\tm}(\beta))+\arg \prod\limits_{\tau=0}^{[\tk s]-1} c_1(\tau/s,\tm,s)=\\=
s\cdot P_{[\tk s]/s,\tm}(\Phi_{\tk,\tm}(\beta))+O(1/s)+\arg \prod\limits_{\tau=0}^{[\tk s]-1} c_1(\tau/s,\tm,s)=\\=
s\cdot P_{\tk,\tm}(\Phi_{\tk,\tm}(\beta))-\{\tk s\}\cdot\left.\frac{\partial P_{\tk_2,\tm}\left(\Phi_{\tk,\tm}(\beta)\right)}{\partial \tk_2}\right|_{\tk_2=\tk}+O(1/s)+\\+s\cdot b_4(\tk,\tm) + b_5(\tk,\tm)+\{\tk s\}\cdot b_6(\tk,\tm).
\end{multline*}
So put 
$$
f^{(s)}_0(\tk,\beta,\tm) := b_5(\tk,\tm)+\{\tk s\}\cdot b_6(\tk,\tm)-\{\tk s\}\cdot \left.\frac{\partial P_{\tk_2,\tm}\left(\Phi_{\tk,\tm}(\beta)\right)}{\partial \tk_2}\right|_{\tk_2=\tk}.
$$
This real scalar is of the form from the statement of our Lemma and enjoys 
\begin{equation*}
	\arg \omega_{\tk,m,s}\left(\Phi_{\tk,\tm}(\beta)\right)=
	\arg w^{{\rm I},s}_{0,\tm}(\beta) +f^{(s)}_0(\tk,\beta,\tm)+ O(1/s).
\end{equation*}
By (\ref{eq:WKB_main}), we have 
$$
	\left|w^{{\rm I},s}_{0,\tm}(\beta)\right| = \frac{1}{\sqrt[4]{Q_{0,\tm}(\beta)}} \cdot\left(1+O(1/s)\right)
$$
and
$$
	\left|\omega_{\tk,m,s}\left(\Phi_{\tk,\tm}(\beta)\right)\right| =e^{b_7(\tk,\tm)+O(1/s)}\cdot \left|w^{{\rm I},s}_{[\tk s]/s,\tm}\left(\Phi_{\tk,\tm}(\beta)\right)\right| 
	 =\frac{e^{b_7(\tk,\tm)}\cdot\left(1+O(1/s)\right)}{\sqrt[4]{Q_{\tk,\tm}(\Phi_{\tk,\tm}(\beta))}}.
$$
Thus, it remains to denote 
$$
f_3(\tk,\beta,\tm) := b_7(\tk,\tm)+\frac{\ln Q_{0,\tm}(\beta) -\ln Q_{\tk,\tm}\left(\Phi_{\tk,\tm}(\beta)\right)}{4},
$$
this quantity obviously satisfies the smoothness condition from the statement because $Q$ does not vanish whenever $|\tm|\le 1/2$. $\blacksquare$

\subsection{Semiclassical measure transformation on cylinder}

This Subsection is devoted to the rigorous proof of heuristic observation from Subsection \ref{subsec:experiment} concerning transformation of semiclassical measure under evolution making waves running with finite velocity. The goal is to prove Proposition \ref{prop:cyl_mapping_exists} stated at the beginning of this Section.

For this, we are going to test semiclassical measure of sequence $\{(u_0^{(s)},1/s)\}_{s\in J}$ by observable $a_0(\beta,\sigma,\eta)\cdot\varphi_4(\xi)$ defined at Subsection \ref{subsubsec:setiing_the_observables} below and show that 
$$
\left\langle (\Op_{1/s} a_0(\beta,\sigma,\eta)\cdot\varphi_4(\xi)) u_0^{(s)}, u_0^{(s)}\right\rangle\approx
\left\langle (\Op_{1/s} a_1(\beta,\sigma,\eta)\cdot\varphi_4(\xi)) u_\tk^{(s)}, u_\tk^{(s)}\right\rangle,
$$
here $u_\tk^{(s)} = \Maas_{0\to \tk}^{(s)} u_0^{(s)}$, whereas  $a_1$ is a special observable defined by (\ref{eq:a1_define}) below and related to $a_0$ by change of variable and multiplication by a given factor.

If someone would find a quantum Hamiltonian $\hat{H}$ such that at least 
$$
\exp\left(\dfrac{\hat{H}\tk}{i\hbar}\right) u_0 \approx \Maas_{0\to \tk}^{(s)} u_0,
$$
say, up to $O(\hbar)$ errors ($\hbar$ is, as usually, something about $1/s$), then the required approximate identity will follow from Yu. Egorov Theorem. But even in this case we stress that we got rid of this Theorem and of Schr\"odinger exponential. This is because we do not have any a priori canonical transformation  derived from Maa\ss{} raising operators and ascension evolution. Mapping $G$ of phase space arises in our study a posteriori as a result of chain of calculations (see \ref{eq:G_def} below). In fact, our proof is just opening the brackets in quadratic form and changing variable in each valuable term.

\subsubsection{Setting the PDO calculus}

\label{subsubsec:quantization}

Let us begin with organizing our quantization procedures.

One starts with quantization in $\R^d, d\in \mathbb N$. A symbol $a=a(x,\xi)$, $x,\xi\in\R^d$, is said to belong to class $\SSS=\SSS(T^*\R^d)$ if any of its partial derivatives is bounded. Take $\hbar>0$ understood as \emph{Planck constant}. For $u\in C_0^\infty(\R^d)$ define \emph{standard quantization} as
\begin{equation}
\label{eq:standard_q_Rd}
(\Op^{\R^d}_\hbar a) u(x) := \int_{\R^d} a(x,\hbar\xi) e^{i\xi\cdot x} \hat u(\xi)\,d\xi,
\end{equation}
where $\hat u$ is the usual Fourier transform. 

Next, we define the action of quantized symbol on  functions $u\in L^2_{\loc}(\R^d)$ without dealing with distributions; the quantization may change a little bit here. Pick any system of non-negative functions $\psi_j\in C_0^\infty(\R^d)$, $j\in \mathbb N$, such that $\sum\limits_{j\in \mathbb N}\psi_j^2\equiv 1$ on $\R^d$ and $\{\supp\psi_j\}_{j\in \mathbb N}$ is locally finite covering of $\R^d$; then quantize by putting 
\begin{equation}
\label{eq:quantize_charts}
(\Op_\hbar a) u := \sum\limits_{j\in \mathbb N}\psi_j\cdot\left(\Op^{\R^d}_\hbar a\right) (\psi_j u) 
\end{equation}
for any $u\in L^2_{\loc}(\R^d)$. \emph{For $a\in \SSS$, the resulting operator does not depend on system $\{\psi_j\}_{j\in \mathbb N}$ up to $O_{L^2_{\loc}(\R^d)\to L^2_{\loc}(\R^d)}(\hbar)$ errors}. This follows from the commutativity of $\SSS\mbox{-operators}$ in the first order and the fact that multiplication by $\psi_j$ is an $\SSS$-operator. Also, the result will change by $O_{L^2_{\loc}(\R^d)\to L^2_{\loc}(\R^d)}(\hbar)$ if we replace "standard" quantization in~(\ref{eq:standard_q_Rd}) by \emph{Weyl} quantization (see \cite{Zw} for details on this way to set the operator calculus).

As above, $\Cyl_l$ is  cylinder with neck of length $l>0$,  and $(\beta,\sigma,\xi,\eta)$ is canonical coordinate system on $T^*\Cyl_l$ (we have $\sigma=\sigma+l$).
We rely to the local charts on $\Cyl_l$ in coordinates $(\beta,\sigma)$; then it makes sense to speak about $\SSS(\Cyl_l)$ and about quantizing symbols from this class by using (\ref{eq:quantize_charts}) with some $\psi_j$'s. This quantization also does not depend on partition of unity up to  $O_{L^2_{\loc}(\Cyl_l)\to L^2_{\loc}(\Cyl_l)}(\hbar)$ errors.

To calculate semiclassical measure on hyperbolic surface $X$, we start with Kohn--Nirenberg symbols of order $0$, the space of such symbols is denoted by $\mathcal S^0(T^*X)$; the (standard or Weyl) quantizations of such symbols are defined correctly as operators up to  $O_{L^2(X)\to L^2(X)}(\hbar)$ freedom. (The advantage of Kohn--Nirenberg symbols is that such a procedure does not depend on local charts and on partition of unity up to operator $O(\hbar)$ errors.) If there exists a closed hyperbolic geodesic of length $l>0$ on $X$ then we cover $X$ by $\Cyl_l$.  Kohn--Nirenberg symbols on $T^*X$ are transferred to $T^*\Cyl_l$ and we may treat quantization on $X$ as quantization on $\Cyl_l$. But then we may decompose symbols from $\mathcal S^0(T^*\Cyl_l)$ into products of symbols from $\SSS(\Cyl_l)$ and this will not lead us to any ambiguity if we restrict ourselves to quantizing in the fixed coordinate system $(\beta,\sigma)$ on $\Cyl_l$.

Let $u\colon\Cyl_l\to\mathbb C$ be a smooth function. Expand it as 
$$
u(\beta,\sigma) = \sum\limits_{m\in 2\pi\mathbb Z/l} \alpha_m(\beta) e^{im\sigma}.
$$
For a symbol  $a\in \SSS(T^*\Cyl_l)$ of the form  $a=a(\beta,\sigma,\eta)$ define another quantization by putting
$$
(\Op^{\mathbb T}_{\hbar} a) u := \sum\limits_{m\in 2\pi\mathbb Z/l} a(\beta,\sigma,m\hbar) \cdot\alpha_m(\beta) e^{im\sigma}.
$$
If $\hbar$ is fixed, this operator depends only on   values of $a$ on a distinguished set of points ($\eta\in 2\pi\hbar\mathbb Z/l$). The philosophy of quantization says that all the reasonable quantizations differ one from another by $O(\hbar)$, and, indeed, we have the following

\begin{lemma}
	\label{lemma:circle_quantize}
	$\Op^{\mathbb T}_{\hbar} a-\Op_{\hbar} a = O_{L^2_{\loc}(\Cyl_l)\to L^2_{\loc}(\Cyl_l)}(\hbar)$. In other words,\vspace{0.051cm} we may apply multiplier to Fourier coefficients instead of Fourier transform of localized function.
\end{lemma}

\noindent A rigorous proof can be given with the help of Paley--Wiener functions used to relate Fourier multipliers on $\R$ and on $\mathbb T$. We omit this step in our exposition.

\subsubsection{Setting the observables}
\label{subsubsec:setiing_the_observables}
Now, we pass to the analysis of eigenfunctions and testing their semiclassical measures. Everywhere we take function $u_0=u_0^{(s)}\colon\Cyl_l\to\R$ such that 
\begin{equation}
\label{eq:eigen2}
-\Delta_{\Cyl_l} u_0 = s^2 u_0.
\end{equation}
Let $\bar\mu^0$ be semiclassical  measure for some subsequence in $\{(u_0^{(s)}, 1/s)\}_{s\in\sqrt{\spec(-\Delta_X)}}$. By Lemma~\ref{lemma:Wigner_usual}, measure $\bar\mu^0$ is supported by the set $\{\cos^2\beta\cdot\left(\xi^2+\eta^2\right)=1\}$. Pick any point $(\beta_0,\sigma_0,\xi_0,\eta_0)\in\Omega^0$ (see (\ref{eq:Omega})).

Fix any function $\varphi_0\in C^\infty_0(\R)$ supported by $[-1/2,1/2]$ such that $0\le\varphi_0\le 1$ on $\R$ and $\varphi_0\equiv 1$ near $0$. Let $\eps>0$ be small enough such that $[\eta_0-\eps,\eta_0+\eps]\subset (-\emm,\emm)$, $[\beta_0-\eps,\beta_0+\eps]\subset (-\pi/2,\pi/2)$. On $T^*\Cyl_l$, define functions 
\begin{gather*}
\varphi_1 = \varphi_1(\eta) := \varphi_{0}\left(\dfrac{\eta-\eta_0}{\eps}\right),\\
\varphi_2 = \varphi_2(\beta) := \varphi_{0}\left(\dfrac{\beta-\beta_0}{\eps}\right),\\
\varphi_3 = \varphi_3(\sigma) := \varphi_{0}\left(\dfrac{\sigma-\sigma_0}{\eps}\right).
\end{gather*}
Put $a_0(\beta,\sigma,\xi,\eta) := \varphi_1(\eta)\varphi_2(\beta)\varphi_3(\sigma)$, this is classical observable from $\SSS(T^*\Cyl_l)$. Pick smooth non-negative $\varphi_4=\varphi_4(\xi)\colon T^*\Cyl_l\to\R$ depending only on $\xi$ which is equal to $1$ on~$\Omega^0$ and for $\xi$ large enough, and equal to $0$ when $\xi<0$; this is possible since $\xi$ is separated from zero on $\Omega_0$, see (\ref{eq:Omega}).  
We are going to express 
$$
\lim\limits_{s\to\infty} \left\langle (\Op_{1/s} a_0\cdot\varphi_4(\xi)) u_{0}, u_{0}\right\rangle
$$
via something similar about $u_\tk$.

\subsubsection{Expansion into $w^{\rm I,\rm II}$-eigenwaves}
\label{subsubsec:cutoff}
Any eigenfunction $u_0$ satisfying (\ref{eq:eigen2}) can be expanded as 
\begin{equation}
\label{eq:eigenfunction_wI_wII_expand}
u_0(\beta,\sigma) = \sum_{m\in 2\pi\mathbb Z/l}e^{im\sigma}\cdot\left(\alpha_m w^{{\rm I},s}_{0, m/s}(\beta)+\alpha^{\rm II}_m w^{{\rm II},s}_{0, m/s}(\beta)\right)
\end{equation}
with some scalars $\alpha_m, \alpha_m^{\rm II}\in\mathbb C$. (Roman superscript in $\alpha_m$ is dropped intentionally to simplify further notation: we will mainly deal with span of $w^{\rm I}$-eigenwaves.) Assume also that $u_0$ is bounded in $L^2_{\loc}(\Cyl_l)$ uniformly by $s$, for this it is enough to suppose that $u_{0}$ projects to a single-valued function on $X$ by covering projection $\pr_{\Cyl_l\to X}$ and that $\int_X |u_0|^2\,d\mathcal A = 1$. 

\begin{lemma} 
\label{lemma:coeff_l2}	
We have
$$
\sum_{\mathclap{\substack{m\in 2\pi\mathbb Z/l\\|m|\le s/2}}} \left(|\alpha_m|^2+|\alpha^{\rm II}_m|^2\right)=O(1)
$$
uniformly by $s$.
\end{lemma}
 
\noindent {\bf Proof.} When $|m|\le s/2$, harmonics $w_{0,\tm}^{j,s}(\beta)$, $j=\rm I, \rm II$, are separated from zero uniformly when $\beta$ ranges any compact subinterval in $\left(-\dfrac{\pi}{2},\dfrac{\pi}{2}\right)$. (Recall that $\tm=m/s$.) So, to prove Lemma, we notice that $u_0\in L^2_{\loc}(\Cyl_l)$ uniformly by $s$, harmonic $e^{im\sigma}\cdot w^{{\rm I},s}_{0,m/s}$ (or $e^{im\sigma}\cdot w^{{\rm II},s}_{0,m/s}$) is orthogonal to $e^{im'\sigma}\cdot w^{{\rm I},s}_{0, m'/s}$ (or $e^{im'\sigma}\cdot w^{{\rm II},s}_{0, m'/s}$) in the direction of any hypercycle $\{\beta\equiv\const\}$ whenever $m\neq m'$ and, finally, $e^{im\sigma}\cdot w^{{\rm I},s}_{0,m/s}$ and $e^{im\sigma}\cdot w^{{\rm II},s}_{0,m/s}$ are almost orthogonal on any hyperbolic geodesic segment of the form  $\{(\beta,\sigma_1)\mid \beta\in (\beta_1,\beta_2)\}$ for some $\beta_1,\beta_2\in (-\pi/2,\pi/2)$, $\sigma_1\in\R$. (The last fact is proved by integrating WKB ansatz (\ref{eq:WKB_main}) for $w$-harmonics by parts.)
$\blacksquare$

\medskip

\noindent We want to cut frequencies with $|m|>s/2$. Take any smooth $\varphi_5 = \varphi_5(\eta)$ supported by $(-1/2,1/2)$ and equal to $1$ near $\supp\varphi_1$. Put $\mathcal P_{1/2} := \Op^{\mathbb T}_{1/s}\varphi_5(\eta)$. One has
$$
\mathcal P_{1/2} u_0 = \sum_{m\in 2\pi\mathbb Z/l}\varphi_5(m/s)e^{im\sigma}\cdot\left(\alpha_m w^{{\rm I},s}_{0, m/s}(\beta)+\alpha^{\rm II}_m w^{{\rm II},s}_{0, m/s}(\beta)\right). 
$$
By the choice of $\varphi_5$ and due to first-order commutativity of quantizations,
$$
\left\langle \Op_{1/s} (a_0\cdot\varphi_4(\xi)) u_{0}, u_{0}\right\rangle =  \left\langle \Op_{1/s} (a_0\cdot\varphi_4(\xi))\circ \mathcal P_{1/2} u_{0}, \mathcal P_{1/2}u_{0}\right\rangle + O(1/s).
$$
So we replace $u_0$ by $\mathcal P_{1/2} u_0$ in the scalar product in the left-hand side of the latter relation.

Let us consider further cut-off, the projection on span of $w^{\rm I}$-waves: put
$$
\mathcal P^{{\rm I},s} \mathcal P_{1/2} u_0 := \sum_{m\in 2\pi\mathbb Z/l}\varphi_5(m/s) \alpha_m \cdot e^{im\sigma}w^{{\rm I},s}_{m/s,0}(\beta).
$$

\begin{lemma}
	\label{lemma:phi4_cut}
	Let $u_0$ be an $s$-eigenfunction bounded in $L^2_{\loc}(\Cyl_l)$ uniformly by $s$.
	\begin{enumerate} 
	\item We have 
	$$
	\Op_{1/s} \varphi_4(\xi)
	\left(\mathcal P_{1/2}u_0-\mathcal P^{{\rm I},s}\mathcal P_{1/2} u_0\right) = O_{L^2_{\loc}(\Cyl_l)}(1/s).
	$$

	\item 
	If $\varphi_4'=\varphi_4'(\xi)\colon \R\to \mathbb C$ is a smooth function compactly supported in $(-\infty,0)$ then 
	$$
	\Op_{1/s} \varphi_4'(\xi) 
	\, \mathcal P^{{\rm I},s}\mathcal P_{1/2} u_0 = O_{L^2_{\loc}(\Cyl_l)}(1/s).
	$$
	
	\item  Semiclassical measure of sequence $\left\{(\mathcal P^{{\rm I},s}\mathcal P_{1/2} u_0^{(s)}, 1/s)\right\}_{s\in J}$ is supported by the set $\{(\beta,\sigma,\xi,\eta)\colon \xi\ge 0\}$.
	\end{enumerate}
\end{lemma}

\noindent The first assertion is intuitively obvious and at least natural: $w^{\rm II}$-functions have negative frequencies by $\xi$ from the explicit WKB ansatz (\ref{eq:WKB_main}), whereas Fourier multiplier $\Op_{1/s} \varphi_4(\xi)$ reserves only positive $\xi$-frequencies. A rigorous proof can be given by inserting $\psi_j = \psi_j^{(1)}(\beta)\cdot\psi_j^{(2)}(\sigma)$ in (\ref{eq:quantize_charts}) with appropriate one-dimensional partitions of unity and then applying Van Der Corput Lemma (or just integrating by parts and noting that phase derivatives are separated from zero). The second assertion of Lemma is analogous to the first one, and the third follows from the second.

Form Lemma \ref{lemma:phi4_cut}, we derive 

\begin{sled}
	\label{sled:phi4_quad_form} Under conditions on $u_0=u_0^{(s)}$ as in Lemma \ref{lemma:phi4_cut}, 
	$$
	\lim\limits_{s\to\infty} \left\langle (\Op_{1/s} a_0\cdot\varphi_4(\xi)) u_0^{(s)}, u_0^{(s)}\right\rangle=
	\lim\limits_{s\to\infty} \left\langle (\Op_{1/s} a_0) \mathcal P^{{\rm I},s} \mathcal P_{1/2}u_0^{(s)}, \mathcal P^{{\rm I},s}\mathcal P_{1/2} u_0^{(s)}\right\rangle.	
	$$
\end{sled}

\noindent 
We replace $\Op$ by $\Op^{\mathbb T}$ using Lemma \ref{lemma:circle_quantize} and transform quadratic form  under limit as 
\begin{equation}
\label{eq:double_sum_init0}
\mathop{\sum\sum}\limits_{m,m'\in 2\pi\mathbb Z/l~} ~ \int\limits_{-\pi/2}^{\pi/2}d\beta\int\limits_{-\infty}^{+\infty} d\sigma\,\alpha_m \bar \alpha_{m'}\,\varphi_1(m/s)\,\varphi_2(\beta)\,\varphi_3(\sigma)\varphi_5(m'/s)\,e^{i(m-m')\sigma}\,w^{{\rm I},s}_{0,m/s}(\beta)\,\bar w^{{\rm I},s}_{0,m'/s}(\beta),
\end{equation}
up to $O(1/s)$ corrections; we used that $\varphi_1(\tm)\cdot\varphi_5(\tm)=\varphi_1(\tm)$ by the choice of $\varphi_5$. (Recall that $\alpha_m$, $m\in 2\pi\mathbb Z/l$, are coefficients from (\ref{eq:eigenfunction_wI_wII_expand})).

\subsubsection{Faraway frequencies}
\label{subsubsec:faraway}
We subdivide double sum (\ref{eq:double_sum_init0}) in two: in the first one frequencies $m,m'$ are far and at the second one they are rather close. The second double sum is more difficult to treat. In this Subsection we show that the first sum --- over faraway frequencies --- can be made negligible.

\begin{lemma}
	\label{lemma:faraway_freqs}
Double sum
	\begin{equation*}
	\mathop{\sum\sum}_{\mathclap{\substack{ m,m'\in 2\pi\mathbb Z/l\\ 
	|m-m'|> s^{1/8}}}} ~~~ \int\limits_{-\pi/2}^{\pi/2}d\beta\int\limits_{-\infty}^{+\infty} d\sigma\,\alpha_m \bar \alpha_{m'}\,\varphi_1(m/s)\,\varphi_2(\beta)\,\varphi_3(\sigma)\,\varphi_5(m'/s)\,e^{i(m-m')\sigma}\,w^{{\rm I},s}_{0,m/s}(\beta)\,\bar w^{{\rm I},s}_{0,m'/s}(\beta) 
	\end{equation*}
	is $O(s^{-1/8})$. 
\end{lemma}

\noindent {\bf Proof.} Integrate by parts by $\sigma$ twice and then apply Young inequality on convolution (or simple Schur test) together with Lemma \ref{lemma:coeff_l2}. $\blacksquare$

\subsection{Nearby frequencies: a version of Yu. Egorov Theorem}

By Lemma \ref{lemma:faraway_freqs}, 
$\left\langle (\Op_{1/s} a_0\cdot\varphi_4(\xi)) u_{0}, u_{0}\right\rangle$
is 
\begin{equation*}
\mathop{\sum\sum}_{\mathclap{\substack{ m,m'\in 2\pi\mathbb Z/l\\ 
			|m-m'|\le s^{1/8}}}} ~~~ \int\limits_{-\pi/2}^{\pi/2}d\beta\int\limits_{-\infty}^{+\infty} d\sigma\,\alpha_m \bar \alpha_{m'}\,\varphi_1(m/s)\,\varphi_2(\beta)\,\varphi_3(\sigma)\,\varphi_5(m'/s)\,e^{i(m-m')\sigma}\,w^{{\rm I},s}_{0,m/s}(\beta)\,\bar w^{{\rm I},s}_{0,m'/s}(\beta) 
\end{equation*}
up to $O(s^{-1/8})$.
For $s$ large, we may drop $\varphi_5$ now and write 
\begin{equation}
\label{eq:double_sum_close_freq}
\mathop{\sum\sum}_{\mathclap{\substack{ m,m'\in 2\pi\mathbb Z/l\\|m|,|m'|\le s/2\\ 
			|m-m'|\le s^{1/8}}}} ~~~ \int\limits_{-\pi/2}^{\pi/2}d\beta\int\limits_{-\infty}^{+\infty} d\sigma\,\alpha_m \bar \alpha_{m'}\,\varphi_1(m/s)\,\varphi_2(\beta)\,\varphi_3(\sigma)\,e^{i(m-m')\sigma}\,w^{{\rm I},s}_{0,m/s}(\beta)\,\bar w^{{\rm I},s}_{0,m'/s}(\beta) 
\end{equation}
instead. 

In this sum, only most valuable terms obtained by opening the brackets remained alive. Now we proceed by changing variables in each term in order to obtain an analogous sum for $u_\tk^{(s)}$. Everywhere we assume that $|m-m'|\le s^{1/8}$ and hence $|\tm-\tm'| = O(s^{-7/8})$ where $\tm=m/s$, $\tm'=m'/s$. 
The assumption on closeness of frequencies will greatly improve our calculations and that is why we dropped terms with faraway frequencies in Lemma \ref{lemma:faraway_freqs}.

\subsubsection{Changing variable in observable}

Write sum (\ref{eq:double_sum_close_freq}) as $\displaystyle\mathop{\sum\sum}_{\mathclap{\substack{ m,m'\in 2\pi\mathbb Z/l\\ |m|, |m'| \le s/2\\ |m-m'|\le s^{1/8}}}} \alpha_m \bar \alpha_{m'} \cdot q_{m,m'}$ with 
\begin{equation}
\label{eq:q0_def}
q_{m,m'} :=  \int\limits_{-\pi/2}^{\pi/2}d\beta\int\limits_{-\infty}^{+\infty} d\sigma\,\varphi_1(m/s)\,\varphi_2(\beta)\,\varphi_3(\sigma)\,e^{i(m-m')\sigma}\,w^{{\rm I},s}_{0,m/s}(\beta)\,\bar w^{{\rm I},s}_{0,m'/s}(\beta). 
\end{equation}

\begin{predl}
	\label{prop:o_1s_correction}
	If we perturb each $q_{m,m'}$ by $O(s^{-3/4})$ then the sum \emph{(\ref{eq:double_sum_close_freq})} will be perturbed by $O(s^{-5/8})$. 
\end{predl}

\noindent {\bf Proof.} By Cauchy--Bunyakovskiy--Schwartz inequality. $\blacksquare$

\medskip

\noindent 
Let's start transforming $q_{m,m'}$.
Introduce new variable $\beta'=\Phi_{\tk,\tm}(\beta)$. Recall that 
\begin{equation*}
\omega_{\tk,m,s}(\beta') = w^{{\rm I},s}_{[\tk s]/s,m/s}(\beta')\cdot \prod\limits_{\tau=0}^{[\tk s]-1} c_1(\tau/s,\tm,s),
\end{equation*}
and, by Lemma \ref{lemma:Phi_diffeo}, 
\begin{equation*}
\omega_{\tk,m,s}\left(\Phi_{\tk,\tm}(\beta)\right)=
w^{{\rm I},s}_{0,\tm}(\beta)\cdot \exp({i f_0^{(s)}(\tk,\beta,\tm)+f_3(\tk,\beta,\tm)+O(1/s)}) 
\end{equation*}
whereas $\Phi_{\tk,\tm}\colon \left(-\dfrac{\pi}2,\dfrac{\pi}2\right) \to  \left(-\dfrac{\pi}2,\dfrac{\pi}2\right)$ is an increasing diffeomorphism depending on $\tm$ smoothly. 
We have
\begin{gather*}
w^{{\rm I},s}_{0,\tm}(\beta)= \omega_{\tk,m,s}\left(\Phi_{\tk,\tm}(\beta)\right)
\cdot \exp\left(-{i f_0^{(s)}(\tk,\beta,\tm)-f_3(\tk,\beta,\tm)+O(1/s)}\right),\\
w^{{\rm I},s}_{0,\tm'}(\beta)= \omega_{\tk,m',s}\left(\Phi_{\tk,\tm'}(\beta)\right)
\cdot \exp\left(-{i f_0^{(s)}(\tk,\beta,\tm')-f_3(\tk,\beta,\tm')+O(1/s)}\right). 
\end{gather*}
Therefore,
\begin{equation*}
\label{eq:w-w_product}
w^{{\rm I},s}_{0,\tm}(\beta) \bar w^{{\rm I},s}_{0,\tm'}(\beta) = \omega_{\tk,m,s}\left(\Phi_{\tk,\tm}(\beta)\right) \cdot \bar  \omega_{\tk,m',s}\left(\Phi_{\tk,\tm'}(\beta)\right) \cdot \exp\left(-2 f_3(\tk,\beta,\tm)+O(s^{-7/8})\right)
\end{equation*}
--- since $|m-m'|\le s^{1/8}$, $|\tm-\tm'| = O(s^{-7/8})$, and $f_0^{(s)}, f_3$ depend on $\tm$ smoothly and take real values, see Lemma \ref{lemma:Phi_diffeo}.
 
We are going to insert this product into (\ref{eq:q0_def}) and change variable as $\beta'=\Phi_{\tk,\tm}(\beta)$. To this end, we need to have $\Phi_{\tk,\tm}(\beta)$ in $\omega_{\tk,m',s}\left(\Phi_{\tk,\tm'}(\beta)\right)$ instead of $\Phi_{\tk,\tm'}(\beta)$. Let's achieve it by calculation of phase correction; \emph{this will shift observable in $\sigma$-direction}.

\begin{predl}
	\label{prop:phase_correction}
 	We have 
 	\begin{equation*}
	\omega_{\tk, m',s}\left(\Phi_{\tk,\tm'}(\beta)\right) = \omega_{\tk, m',s}\left(\Phi_{\tk,\tm}(\beta)\right) \cdot \exp\left(i f_4(\tk,\beta,\tm) (m'-m)+O(s^{-3/4})\right)
	\end{equation*}
with
$$
f_4(\tk,\beta,\tm)=
\left(\tk+\sqrt{Q_{\tk,\tm}\left(\Phi_{\tk,\tm}(\beta)\right)}\right)\cdot\frac{\partial \Phi_{\tk,\tm}(\beta)}{\partial\tm}.
$$	
\end{predl}

\noindent {\bf Proof.}
By (\ref{eq:omega_def}),
\begin{equation*}
\omega_{\tk,m',s}\left(\Phi_{\tk,\tm'}(\beta)\right) = w_{[\tk s]/s,\tm'}^{{\rm I},s}\left(\Phi_{\tk,\tm'}(\beta)\right)\cdot \prod_{\tau=0}^{[\tk s]-1} c_1(\tau/s,\tm,s).
\end{equation*}
WKB ansatz (\ref{eq:WKB_main}) for the first multiplier is
\begin{multline*}
 w_{[\tk s]/s,\tm'}^{{\rm I},s}\left(\Phi_{\tk,\tm'}(\beta)\right) = \\=
 \exp\left(is\cdot\left([\tk s]/s\cdot\Phi_{\tk,\tm'}(\beta)+ \int\limits_0^{\Phi_{\tk,\tm'}(\beta)} \sqrt{Q_{[\tk s]/s,\tm'}}\right) -\frac14\int\limits_0^{\Phi_{\tk,\tm'}(\beta)}\frac{Q_{[\tk s]/s,\tm'}'}{Q_{[\tk s]/s,\tm'}}+O\left(s^{-1}\right)\right).
\end{multline*}
By Taylor formula, the latter is 
\begin{multline*}
w_{[\tk s]/s,\tm'}^{{\rm I},s}\left(\Phi_{\tk,\tm'}(\beta)\right) =\\=  
w_{[\tk s]/s,\tm'}^{{\rm I},s}\left(\Phi_{\tk,\tm}(\beta)\right) \cdot \exp\left(i(m'-m)\left(\tk+\sqrt{Q_{\tk,\tm}\left(\Phi_{\tk,\tm}(\beta)\right)}\right)\cdot\frac{\partial \Phi_{\tk,\tm}(\beta)}{\partial\tm}+O(s^{-3/4})\right).
\end{multline*}
This leads to the desired. $\blacksquare$

\medskip

\noindent So, we have
\begin{multline*}
w^{{\rm I},s}_{0,\tm}(\beta) \bar w^{{\rm I},s}_{0,\tm'}(\beta) =\\= \omega_{\tk,m,s}\left(\Phi_{\tk,\tm}(\beta)\right) \cdot \bar  \omega_{\tk,m',s}\left(\Phi_{\tk,\tm}(\beta)\right)\cdot \exp\left(-2 f_3(\tk,\beta,\tm)+i(m-m') f_4(\tk,\beta,\tm)+O(s^{-3/4})\right).
\end{multline*}
Inserting this into (\ref{eq:q0_def}) and changing variable as $\beta'=\Phi_{\tk,\tm}(\beta)$, $\beta=\Phi_{\tk,\tm}^{-1}(\beta')$, we get
\begin{multline*}
q_{m,m'} =
O(s^{-3/4})+\int\limits_{-\pi/2}^{\pi/2}d\beta'\int\limits_{-\infty}^{+\infty} d\sigma\,\frac{\partial\Phi_{\tk,\tm}^{-1}(\beta')}{\partial\beta'}\cdot\varphi_1(\tm)\, \varphi_2\left(\Phi_{\tk,\tm}^{-1}(\beta')\right)\varphi_3(\sigma)\,\times\\\times\exp\left(-2 f_3(\tk,\Phi_{\tk,\tm}^{-1}(\beta'),\tm)\right)\cdot\exp\left[i(m-m')(\sigma+f_4(\tk,\Phi_{\tk,\tm}^{-1}(\beta'),\tm))\right]\omega_{\tk,m,s}(\beta') \cdot \bar  \omega_{\tk,m',s}(\beta')=\\=
O(s^{-3/4})+\int\limits_{-\pi/2}^{\pi/2}d\beta'\int\limits_{-\infty}^{+\infty} d\sigma\,\frac{\partial\Phi_{\tk,\tm}^{-1}(\beta')}{\partial\beta'}\cdot\varphi_1(\tm)\, \varphi_2\left(\Phi_{\tk,\tm}^{-1}(\beta')\right)\varphi_3\left(\sigma-f_4(\tk,\Phi_{\tk,\tm}^{-1}(\beta'),\tm)\right)\,\times\\\times \exp\left(-2 f_3(\tk,\Phi_{\tk,\tm}^{-1}(\beta'),\tm)\right)\cdot \exp\left(i(m-m')\sigma\right)\cdot \omega_{\tk,m,s}(\beta') \cdot \bar  \omega_{\tk,m',s}(\beta').
\end{multline*}
Thus, if we put
\begin{multline}
\label{eq:a1_define}
a_1(\beta',\sigma,\eta) := \\:=\frac{\partial\Phi_{\tk,\eta}^{-1}(\beta')}{\partial\beta'}\cdot \varphi_1(\eta)\, \varphi_2\left(\Phi_{\tk,\eta}^{-1}(\beta')\right)\varphi_3\left(\sigma-f_4(\tk,\Phi_{\tk,\eta}^{-1}(\beta'),\eta)\right)\cdot \exp{\left(-2 f_3(\tk,\Phi_{\tk,\eta}^{-1}(\beta'),\eta)\right)},
\end{multline}
then, by Proposition \ref{prop:o_1s_correction}, we have

\begin{predl}
	\label{prop:a1_observe}
	Sum \emph{(\ref{eq:double_sum_close_freq})} is 
	$$
	O(s^{-5/8})+\mathop{\sum\sum}_{\mathclap{\substack{ m,m'\in 2\pi\mathbb Z/l\\ |m|, |m'| \le s/2\\ |m-m'|\le s^{1/8}}}} ~~~ \int\limits_{-\pi/2}^{\pi/2}d\beta'\int\limits_{-\infty}^{+\infty} d\sigma\,\alpha_m\, \bar \alpha_{m'}\cdot a_1(\beta',\sigma,m/s) e^{i(m-m')\sigma}	\omega_{\tk,m,s}(\beta') \cdot \bar  \omega_{\tk,m',s}(\beta').
	$$
\end{predl}

\subsubsection{Reverting the observable}

Classical observable $a_1(\beta,\sigma,\eta)$ was defined above in (\ref{eq:a1_define}), and $\varphi_4(\xi)$ was introduced in Subsection \ref{subsubsec:setiing_the_observables}. We have also defined $u_\tk := \Maas_{0\to \tk}^{(s)} u_0$. Consider new quadratic form 
\begin{equation}
\label{eq:q1_init}
\left\langle (\Op_{1/s} a_1\cdot\varphi_4(\xi)) u_\tk, u_\tk\right\rangle
\end{equation}
and transform it to the expression similar to (\ref{eq:double_sum_close_freq}). For this, write $u_\tk$ as
\begin{equation*}
u_\tk(\beta,\sigma) = \sum_{m\in 2\pi\mathbb Z/l}e^{im\sigma}\cdot\left(\alpha_{\tk,m} w^{{\rm I},s}_{[\tk s]/s,m/s}(\beta)+\alpha^{\rm II}_{\tk,m} w^{{\rm II},s}_{[\tk s]/s,m/s}(\beta)\right). 
\end{equation*}
For any fixed $m$, $|m|\le s/2$, by Corollary \ref{sled:exc_coeff}, 
\begin{equation}
\label{eq:alpha_m_1}
\alpha_{\tk,m} = \alpha_m \cdot \prod\limits_{\tau=0}^{[\tk s]-1} c_1(\tau/s,\tm,s) + O(1/s)\cdot \left(|\alpha_m|+|\alpha_m^{\rm II}|\right)
\end{equation}
with $c_1(\tau/s,\tm,s)$ as in Lemma \ref{lemma:monochrome}. 

Arguing as in Subsections \ref{subsubsec:cutoff} and \ref{subsubsec:faraway} above, we transform (\ref{eq:q1_init}) to 
\begin{equation}
\label{eq:u_1_double_sum_close_freq}
\mathop{\sum\sum}_{\mathclap{\substack{ m,m'\in 2\pi\mathbb Z/l\\ |m|, |m'| \le s\emm\\ |m-m'|\le s^{1/8}}}} ~~~ \int\limits_{-\pi/2}^{\pi/2}d\beta\int\limits_{-\infty}^{+\infty} d\sigma\,\alpha_{\tk,m} \bar \alpha_{\tk,m'}\,a_1(\beta,\sigma,m/s)\,e^{i(m-m')\sigma}\,w^{{\rm I},s}_{[\tk s]/s,m/s}(\beta)\,\bar w^{{\rm I},s}_{[\tk s]/s,m'/s}(\beta),
\end{equation}
up to $O(s^{-1/8})$ errors (now we cut all the frequencies outside of $(-\emm,\emm)$, recall that $\supp\varphi_1$ lies in this interval). Indeed,  inserting $\varphi_4$ preserves only $w^{\rm I}$-waves because main frequency term in WKB ansatz (\ref{eq:WKB_main}) for $w^{\rm II}$-waves is $\tk-\sqrt{2\tk\tm\tg\beta-\tm^2+1+\tg^2\beta+\tk^2}$ which is strictly negative and separated from zero, this is true by the choice of $\emm$ in (\ref{eq:eta_max_def}) and by Cauchy--Schwartz inequality. The other steps are also similar to the case of~$u_0$. 

Next, using (\ref{eq:alpha_m_1}), we rewrite (\ref{eq:u_1_double_sum_close_freq}) as 
\begin{equation*}
\mathop{\sum\sum}_{\mathclap{\substack{ m,m'\in 2\pi\mathbb Z/l\\ |m|, |m'| \le \emm\\ |m-m'|\le s^{1/8}}}} ~~~ \int\limits_{-\pi/2}^{\pi/2}d\beta\int\limits_{-\infty}^{+\infty} d\sigma\,\alpha_{m} \bar \alpha_{m'}\,a_1(\beta,\sigma,m/s)\,e^{i(m-m')\sigma}\,\omega_{\tk,m,s}(\beta)\,\bar \omega_{\tk,m',s}(\beta). 
\end{equation*}
up to $O(s^{-7/8})$ correction. But this is the form which we were left with in the previous Subsection (Proposition \ref{prop:a1_observe}). We thus have 
	\begin{equation*}
	\left\langle (\Op_{1/s} a_0(\beta,\sigma,\eta)\cdot\varphi_4(\xi)) u_0^{(s)}, u_0^{(s)}\right\rangle=O(s^{-1/8})+
	\left\langle (\Op_{1/s} a_1(\beta,\sigma,\eta)\cdot\varphi_4(\xi)) u_\tk^{(s)}, u_\tk^{(s)}\right\rangle.
	\end{equation*}
By limit pass,
$$
\int\limits_{T^*\Cyl_l} a_0 \cdot\varphi_4\, d\bar\mu^0=
\int\limits_{T^*\Cyl_l} a_1\cdot\varphi_4 \, d\bar\mu^\tk.
$$
Note that $(\beta,\sigma,\eta)$ is a coordinate system on $\Omega^0$ and on $\Omega^\tk$. In these coordinates, put
\begin{equation}
\label{eq:G_def}
G \colon \Omega^0\to \Omega^{\tk}, ~~~
G(\beta,\sigma,\eta) := \left(\Phi_{\tk,\eta}(\beta),\sigma+f_4(\tk,\beta,\eta),\eta\right), ~~~ (\beta,\sigma,\eta)\in\Omega^0,
\end{equation}
and also define $A\colon \Omega^\tk\to(0,+\infty)$ by 
\begin{equation*}
A(\beta,\sigma,\eta) := \left(\frac{\partial\Phi_{\tk,\eta}^{-1}(\beta)}{\partial\beta}\right)^{-1}\cdot\exp\left(2f_3(\tk,\Phi_{\tk,\eta}^{-1}(\beta),\eta)\right), ~~  (\beta,\sigma,\eta)\in\Omega^\tk.
\end{equation*}
Since $a_0$ was arbitrary observable of the form from Subsection \ref{subsubsec:setiing_the_observables}, and $a_1$ is related to $a_0$ by (\ref{eq:a1_define}), we then have $\bar\mu^{\tk}=A\cdot G_\sharp\bar\mu^0$ up to restrictions to $\Omega^\tk$ and $\Omega^0$ respectively. The proof of Proposition \ref{prop:cyl_mapping_exists} is complete.~$\blacksquare$

\subsection{Testing the mapping}
\label{subsec:testing}

We proceed investigation of semiclassical measure transfer \emph{on cylinder $\Cyl_l$}. To prove Theorem \ref{th:finite_time_shift}, it remains to test mappings $\overline T\mu= A\cdot G_\sharp\mu$, $\mu \in \Meas(\Omega^0)$, and $\underline T := (\phi_{\tk})_\sharp\overline T(\phi_{0}^{-1})_\sharp\colon \Meas((\phi_0)_\sharp\Omega^0)\to\Meas((\phi_\tk)_\sharp\Omega^\tk)$ by substituting there semiclassical measures of appropriately concentrated functions. (Mapping $G$ and function $A$ were defined at the end of the previous Subsection.) If $\mu$ as above is geodesic line then $\overline{T}\mu$ should be a $\tk$-hypercycle with a scalar coefficient. By calculation of asymptotics of $G$, we recover this hypercycle (it is enough to find its ends on the absolute of $\Cyl_l$). To recover the scalar coefficient before $\tk$-hypercycle, we substitute to our result a quantum ergodic sequence which does exist by Shnirel'man--Zelditch--Colin de Verdi\`ere Theorem (\cite{Shn74}, \cite{ZelditchQE}, \cite{CdV}).

First, let us prove that there do exist semiclassical measures concentrated on geodesics. 
Recall that magnetic Hamiltonian on cylinder is 
$$
H_\tk = \frac{(\xi-\tk)^2\cos^2\beta +(\eta\cos\beta-\tk\sin\beta)^2}2,
$$
where $(\beta,\sigma,\xi,\eta)$ is the canonical coordinate system on $\Cyl_l$. Thence, $\eta$ remains constant under Hamiltonian evolution $\exp t\Xi_\tk$ on $T^*\Cyl_l$.

Let $\underline{\gamma}=\underline{\gamma}(t)$, $t\in\R$, be parametrized geodesic on $\Cyl_l$ intersecting neck of cylinder $\{\beta=0\}$ transversally and such that $\beta$ increases along $\underline{\gamma}$; there exists a plenty of such $\underline{\gamma}$, which can be seen by considering geodesics in coordinates $(\beta,\sigma)$ on $\HH\simeq\mathbb C^+$ and then by folding the latter plane to cylinder.  Put $\bar{\gamma}(t):= (\phi^{-1}_0)\underline\gamma'(t)$. We have $\eta=\const=:\eta_0$ along $\bar{\gamma}$; we always assume that $|\eta_0|<\emm$. Also, $\xi>0$ along $\bar{\gamma}$ since $0<\dot{\beta}=\xi\cdot\cos^2\beta$ there. Denote by $\bar{\mu}^{\bar\gamma}$ the positive measure on $T^*\Cyl_l$ supported by $\bar{\gamma}$, invariant with respect to geodesic flow $\exp t\Xi_0$ and normed such that lifts of length $1$ segments on $\underline{\gamma}$ have unit mass.

\begin{lemma}
	\label{lemma:single_geodesic_exist}
	Let $\bar{\gamma}\subset T^*\Cyl_l$ be as above. 
	
	There exists a sequence of positive scalars $\{s\}_{s\in J}$ going to $+\infty$ and a sequence of functions $u_0^{(s)} \colon \Cyl_l\to\R$ bounded in $L^2_{\loc}(\Cyl_l)$ uniformly by $s$ and enjoying $-\Delta_{\Cyl_l} u_0^{(s)}=s^2 u_0^{(s)}$ such that the semiclassical measure of the sequence $\{(u_0^{(s)},1/s)\}_{s\in J}$ is $\bar{\mu}^{\bar{\gamma}}$.
\end{lemma}

\noindent {\bf Proof.} Without loss of generality, we may assume that $\underline{\gamma}$  passes through the point $(\beta,\sigma)=(0,0)$. 
Take $s=1,2,\dots$, and let $\{K(s)\}_{s\in\mathbb N}$ be real sequence increasing to $+\infty$ slowly enough. Define
$$
I(s) := \left\{k\in\R\colon \eta_0s+k\in 2\pi\mathbb Z/l, \, |k|\le \sqrt{K(s)}\right\}.
$$
Put 
$$
u_0^{(s)}(\beta,\sigma) := \frac{1}{\sqrt{K(s)}}\cdot\sum\limits_{k\in I(s)} e^{i(k+\eta_0s)\sigma} \cdot w_{0,\eta_0+k/s}^{{\rm I},s}(\beta), ~~~ (\beta,\sigma)\in\Cyl_l.
$$
Here, $\eta_0$ is the constant value of $\eta$ along $\bar{\gamma}$ whereas $w_{0,\eta_0+k/s}^{{\rm I},s}$ is cylindric harmonic from~(\ref{eq:WKB_main}) with $+$ sign. Then $u_0^{(s)}$ is defined on $\Cyl_l$ as a single-valued function, satisfies $-\Delta_{\Cyl_l}u_0^{(s)}=s^2 u_0^{(s)}$, and also $u_0^{(s)}\in L^2_{\loc}(\Cyl_l)$ uniformly by $s$, this is by WKB ansatz~(\ref{eq:WKB_main}). Also, by (\ref{eq:WKB_main}) we may write
\begin{multline*}
u_0^{(s)}(\beta,\sigma) = \frac{1}{\sqrt{K(s)}}\cdot\sum\limits_{k\in I(s)} \exp\left({i(k+\eta_0s)\sigma+is f_5(\beta,\eta_0)+ikf_6(\beta,\eta_0)+O(s^{-1/2})}\right)=\\=\frac{ e^{is \left(f_5(\beta,\eta_0)+\eta_0\sigma\right)}}{\sqrt{K(s)}}\cdot\sum\limits_{k\in I(s)} \exp\left({ik(\sigma+f_6(\beta,\eta_0))+O(s^{-1/2})}\right), ~~ (\beta,\sigma)\in\Cyl_l,
\end{multline*}
with smooth $f_5,f_6$ satisfying $f_5(0,\eta_0)=f_6(0,\eta_0)=0$; this holds if $K(s)$ does not grow too fast. But then any weak* limit of subsequence in $\{|u_0^{(s)}|^2\}_{s\in\mathbb N}$ has to be concentrated on the set $$\{(\beta,\sigma)\in\Cyl_l\colon  \sigma=-f_6(\beta,\eta_0)\}$$
(outside of this set, function $u_0^{(s)}$ becomes pointwise small uniformly on compacta when $s$ is large).
Also, any semiclassical measure of a subsequence in $\{(u_0^{(s)},1/s)\}_{s\in\mathbb N}$ is $\exp t\Xi_0$-invariant and is supported by the set $\{\eta=\eta_0\}\subset T^*\Cyl_l$. (The latter can be proved by using Lemma \ref{lemma:circle_quantize} and $\Op^{\mathbb T}$-quantization by circular sections of $\Cyl_l$.) Also, this measure is supported by $\{\xi\ge 0\}$, that is because all $w^{\rm I}$-harmonics have positive $\xi$-frequencies. Under all these conditions, such a semiclassical measure can be only $\bar\mu^{\bar\gamma}$ up to a positive scalar factor. We may rescale initial functions $u_0^{(s)}$ to obtain exactly $\bar\mu^{\bar\gamma}$. Proof is complete.~$\blacksquare$

\medskip

\noindent 
Take $\underline{\mu}^{0} := (\phi_0)_\sharp \bar{\mu}^{\bar\gamma}$ with $\bar{\mu}^{\bar{\gamma}}$ as in Lemma \ref{lemma:single_geodesic_exist}. For $u_0^{(s)}$ as in this Lemma, put $u_\tk^{(s)} := \Exc_{0\to \tk}^{(s)} u_0^{(s)}$ and let $\bar{\mu}^{\tk}$ be semiclassical measure of a subsequence in $\{(u_\tk^{(s)},1/s)\}_{s\in J}$. Put also $\underline{\mu}^\tk := (\phi_\tk)_\sharp\bar{\mu}^{\tk}$.

By Proposition \ref{prop:cyl_mapping_exists} which is already proven, measure $\bar{\mu}^{\tk}$ is concentrated on a smooth curve in $T^*\Cyl_l$. Also, $\bar{\mu}^{\tk}$ is $\exp t\Xi_\tk$-invariant. Then $\underline{\mu}^\tk$ is a measure concentrated on some $\tk$-hypercycle $\underline\gamma_\tk$. Let us reconstruct $\underline\gamma_\tk$ by its {ideal points}. Recall that $(\beta,\sigma,\eta)$ is coordinate system on $\Omega^\tk$.

\begin{predl}
\label{prop:bounded_dist}
	Let $\tk,\eta$ be fixed. When $\beta$ tends to $\pm\pi/2$, hyperbolic $\Cyl_l$-distance between $(\beta,\sigma)$ and basepoint of covector $G(\beta,\sigma,\eta)$ 	remains bounded.
\end{predl}

\noindent {\bf Proof.} Let $\beta\to\pi/2$, the other case is similar. By (\ref{eq:G_def}), the basepoint of $G(\beta,\sigma,\eta)$ from the statement has coordinates $\beta'=\Phi_{\tk,\eta}(\beta)$, $\sigma'=\sigma+f_4(\tk,\beta,\eta)$. By the properties of $\Phi$, we have $\beta'\to\pi/2$ when $\beta\to\pi/2$.

First, we deal with $\beta'$. By the proof of Lemma \ref{lemma:Phi_diffeo}, 
\begin{equation}
\label{eq:Phi_def_again}
\tk\beta'+\int\limits_{0}^{\beta'} \sqrt{Q_{\tk,\eta}(\beta_1)}\,d\beta_1 = \int\limits_{0}^{\beta} \sqrt{Q_{0,\eta}(\beta_1)}\,d\beta_1-b_4(\tk,\eta)
\end{equation}
with smooth $b_4$ \emph{not depending on $\beta$} (see (\ref{eq:Phi_define})). Recall that 
$$
Q_{\tk,\eta}(\beta)=2\tk\eta\tg\beta-\eta^2+\dfrac{1}{\cos^2\beta}+\tk^2,
$$
this explodes when $\beta\to\pi/2$. By (\ref{eq:Phi_def_again}), $$\int\limits_{0}^{\beta'} \sqrt{Q_{\tk,\eta}(\beta_1)}\,d\beta_1 - \int\limits_{0}^{\beta} \sqrt{Q_{0,\eta}(\beta_1)}\,d\beta_1$$
remains bounded when $\beta,\beta'\to\pi/2$. Thence  $\ln\left(\frac{\pi}{2}-\beta'\right)- \ln\left(\frac{\pi}{2}-\beta\right)$ remains bounded with $\beta,\beta'\to\pi/2$, that is, $\frac{\pi}{2}-\beta'$ is comparable to $\frac{\pi}{2}-\beta$. By taking integral, the hyperbolic distance between $(\beta,\sigma)$ and $(\beta',\sigma)$ is $\left|\ln\left(\frac{1+\sin\beta}{\cos\beta}\right)-\ln\left(\frac{1+\sin\beta'}{\cos\beta'}\right)\right|$ and stays bounded when $\beta\to\pi/2$.

Now we estimate shift in $\sigma$-direction. For the $f_4$,  Proposition \ref{prop:phase_correction} gives the expression 
$$
f_4(\tk,\beta,\eta)= 
\left(\tk+\sqrt{Q_{\tk,\eta}\left(\Phi_{\tk,\eta}(\beta)\right)}\right)\cdot\frac{\partial \Phi_{\tk,\eta}(\beta)}{\partial\eta}.
$$
The first factor is comparable to $\dfrac{1}{\pi/2-\beta}$. For the second one, we differentiate (\ref{eq:Phi_def_again}) by $\eta$ ($\beta'=\Phi_{\tk,\eta}(\beta)$ therein) and get
\begin{equation}
\label{eq:Phi_diff}
\frac{\partial \Phi_{\tk,\eta}(\beta)}{\partial\eta} = \dfrac
{
	\int\limits_0^{\beta}\dfrac{
		\partial\sqrt{Q_{0,\eta}(\beta_1)}}{\partial\eta}\,d\beta_1-
	\int\limits_0^{\Phi_{\tk,\eta}(\beta)}\dfrac{\partial\sqrt{Q_{\tk,\eta}(\beta_1)}}{\partial\eta}\,d\beta_1
	-\dfrac{\partial b_4(\tk,\eta)}{\partial\eta}
}
{
	\tk+\sqrt{Q_{\tk,\eta}(\Phi_{\tk,\eta}(\beta))}
}.
\end{equation}
Denominator at the right-hand side grows as $\dfrac{1}{\pi/2-\beta}$ with $\beta\to\pi/2$. We claim that numerator is $O(\pi/2-\beta)$ for $\beta$ close to $\pi/2$.

From the proof of Lemma \ref{lemma:Phi_diffeo}, we have $b_4(\tk,\eta) = \int\limits_0^\tk b_1(\tk_2,\eta)\,d\tk_2$, $b_1(\tk_2,\eta)$ is argument of main part of $c_1(\tk_2,\eta,s)$, that is  $-\arctg\frac{\eta}{\sqrt{\tk_2^2-\eta^2+1}}$. By a calculation we have $\dfrac{\partial b_4(\tk,\eta)}{\partial\eta}=\ln\left(\sqrt{1-\eta^2}\right)-\ln\left(\tk+\sqrt{\tk^2-\eta^2+1}\right)$. Derivatives of square roots in (\ref{eq:Phi_diff}) are bounded, we conclude that numerator in (\ref{eq:Phi_diff}) is
\begin{multline*}
\int\limits_0^{\pi/2}\dfrac{\eta\cos\beta_1-\tk\sin\beta_1}{\sqrt{\tk^2\cdot\cos^2\beta_1+1-\eta^2\cdot\cos^2\beta_1+2\tk\eta\cdot\sin\beta_1\cos\beta_1}}\,d\beta_1+\ln\left(\tk+\sqrt{\tk^2-\eta^2+1}\right)-\\-\left(\mbox{the same  at }\tk=0\right)+O(\pi/2-\beta).
\end{multline*}
By substituting $B = \eta\sin\beta_1+\tk\cos\beta_1$ we see that the integral here is $\ln\left(\eta+1\right)-\ln\left(\tk+\sqrt{\tk^2-{\eta}^{2}+1}\right).$ Thus, only $O(\pi/2-\beta)$ remains alive in numerator of right-hand side in (\ref{eq:Phi_diff}).

Gathering all estimates, we see that $f_4(\tk,\beta,\eta)=O\left(\pi/2-\beta\right)$. Since metric tensor on $\Cyl_l$ is $ds^2=\cos^{-2}\beta\cdot (d\beta^2+d\sigma^2)$, we conclude that distance between $(\beta',\sigma)$ and $(\beta',\sigma')=(\beta',\sigma+f_4(\tk,\beta,\eta))$ stays bounded with $\beta\to\pi/2$. The same holds when $\beta\to-\pi/2$. Proof is complete.
$\blacksquare$

\medskip

\noindent So, $\bar\mu^\tk$ has to be concentrated on a $\tk$-hypercycle in $T^*\Cyl_l$ whose projection to $\Cyl_l$ has the same ideal points as $\underline{\gamma}$. On $T^*\Cyl_l$, there exist only two such $\tk$-hypercycles, one of them is   $\phi_\tk^{-1}\mathcal T_{\tk}\underline \gamma' =\phi_\tk^{-1}\circ \Sc_{\sqrt{\tk^2+1}}\circ \mathcal{R}_{\pi/2}\circ h^0_{\ln\left(\tk+\sqrt{\tk^2+1}\right)}\circ\mathcal{R}_{-\pi/2} \gamma'$, and the second is $\phi_\tk^{-1}\mathcal T_{\tk} \mathcal{R}_\pi\underline\gamma'$. 
 The latter is not possible since all our constructions are continuous with respect to $\tk$, and for $\tk=0$ we have identical transformation of measure (which is not $\mathcal R_{\pi}$-rotated).

Let us summarize what we have. Denote by $\underline{\mu}_\tk$ the measure supported by $\mathcal T_{\tk}\underline \gamma'$, invariant with respect to $\tk$-hypercyclic flow $h_\tk$ and normed such that $\tk$-hypercyclic segments of length $1$ have unit mass. We have already proved that 
$$
\underline{T}\underline\mu^{0}=f_7(\eta_0)\cdot \underline{\mu}_\tk
$$
with some smooth positive $f_7$. It remains to check that $f_7$ does not depend on $\eta_0$ ($|\eta_0|<\emm$). 

For this, consider any compact hyperbolic surface $X$ and a quantum ergodic sequence $u_0^{(s)}$ $(s\in J)$ of functions on $X$ with $-\Delta_X u_0^{(s)}= s^2 u_0^{(s)}$; this means that if $\bar{\mu}^0$ is semiclassical measure of sequence $\{(u_0^{(s)},1/s)\}_{s\in J}$ then $(\phi_0)_\sharp\bar{\mu}^0$ is the uniform Liouville measure on $S_1X$ (see definition after Proposition \ref{prop:horocyclic_que}). Such a sequence exists by \cite{Shn74}, \cite{ZelditchQE}, \cite{CdV}. If $\bar{\mu}^\tk$ is semiclassical measure of sequence $\{(\Exc_{0\to \tk}^{(s)} u_0^{(s)},1/s)\}_{s\in J}$ constructed via some covering of $X$ by $\HH$ then $(\phi_\tk)_\sharp \bar{\mu}^\tk$ does not depend on covering (Proposition~\ref{prop:invariance}) and thus can be considered as a measure on $S_{\sqrt{\tk^2+1}}X$. Also, $(\phi_\tk)_\sharp \bar{\mu}^\tk$ is invariant with respect to $\tk$-hypercyclic flow on $S_{\sqrt{\tk^2+1}} X$. By Proposition \ref{prop:cyl_mapping_exists}, $(\phi_\tk)_\sharp \bar{\mu}^\tk$ is absolutely continuous with respect to coordinates in $S_{\sqrt{\tk^2+1}} X$; $\tk$-hypercyclic flow is ergodic on this level set since it is conjugate to the geodesic flow. Thence, $(\phi_\tk)_\sharp \bar{\mu}^\tk$ is constant times the uniform measure on $S_{\sqrt{\tk^2+1}} X$.  If, on $X$, there exists a closed geodesic loop of length $l$ then measures $(\phi_0)_\sharp \bar{\mu}^0$ and $(\phi_\tk)_\sharp \bar{\mu}^\tk$ can be transferred to $T\Cyl_l$. 
Testing transformation $\underline T$ by uniform measures, we conclude that $f_7(\eta_0)$ does not depend on $\eta_0$. Nor it depends on~$l$, the necklength of $\Cyl_l$, since $A$ and $G$ do not,  this can be seen from all the constructions in this Section. 
Then, on cylinder we obtain, now for $\bar{\mu}^0$ and $\bar{\mu}^\tk$ constructed from an arbitrary sequence of eigenfunctions, that 
$$
\left(\Sc_{\sqrt{\tk^2+1}}\circ \mathcal{R}_{\pi/2}\circ h^0_{\ln\left(\tk+\sqrt{\tk^2+1}\right)}\circ\mathcal{R}_{-\pi/2}\right)_\sharp (\mathds 1_{\Omega^0}\cdot\bar{\mu}^0)=C\cdot\mathds 1_{\Omega^\tk}\cdot\bar{\mu}^\tk
$$
with an absolute constant $C$.

Closed geodesics are dense in $S_1X$ (see, e.g. \cite{KatokHasselblat}). Then the union of sets of the form $\pr_{T\Cyl_l\to TX} \phi_\tk\Omega^0$ constructed by all cylindric coverings of $X$ is all the $S_1X$. Then, at $X$, we have
$$
\left(\Sc_{\sqrt{\tk^2+1}}\circ \mathcal{R}_{\pi/2}\circ h^0_{\ln\left(\tk+\sqrt{\tk^2+1}\right)}\circ\mathcal{R}_{-\pi/2}\right)_\sharp (\bar{\mu}^0)=C\cdot\bar{\mu}^\tk.
$$
But $C=1$ since $\int_X|u_\tk|^2\,d\mathcal A$ does not depend on $\tk$. Theorem \ref{th:finite_time_shift} is proved. $\blacksquare$

\begin{sled}
	QUE for functions $u_0^{(s)}$, $s\in\sqrt{\spec(-\Delta_X)}$, is equivalent to QUE for functions $u_\tk^{(s)}$.
\end{sled}

\section{Infinite ascension}

We know that semiclassical measures of functions $\Exc_{0\to \tk}^{(s)} u_0^{(s)}$ are invariant under $\tk\mbox{-hypercyclic}$ flow. This means, roughly speaking, that measures $\left|\Exc_{0\to \tk}^{(s)} u_0^{(s)}\right|^2\,d\mathcal A$ are decomposable into $\tk$-hypercycles, the curves of curvature $\tk/\sqrt{\tk^2+1}$. The latter tends to $1$ when $\tk\to+\infty$. Thus, measures $\left|\Exc_{0\to \tk}^{(s)} u_0^{(s)}\right|^2\,d\mathcal A$ are \emph{almost} decomposable into horocycles when $\tk\to+\infty$. Due to Furstenberg's Theorem on unique ergodicity of horocyclic flow,   we may expect chaotic behavior of $\Exc_{0\to \tk}^{(s)} u_0^{(s)}$ for large $\tk$. So, in this Subsection we prove horocyclic QUE. Intuitively, we have a good chance to succeed: unique ergodicity of classic horocyclic flow must have some quantum counterpart by Bohr correspondence principle.

In fact, we may forget  the origin of functions $\Exc_{0\to \tk}^{(s)} u_0^{(s)}$ and prove the following

\begin{predl}
	\label{prop:horocyclic_que}
	Let $\{s_n\}_{n=1}^\infty\subset \R\setminus\{0\}$, $\{\tau_n\}_{n=1}^\infty\subset \mathbb Z$ be any sequences with the only assumption that 
	$$
	\frac{\tau_n}{s_n}\xrightarrow{n\to\infty}\infty.
	$$
	Let $\Gamma<\Isom^+(\HH)$ be a torsion-free group with a compact fundamental domain $F$. Suppose that functions $U_n\colon\HH\to\mathbb C$, $n=1,2,\dots$, are such that $U_n\in\mathcal F^{\tau_n}(\Gamma)$, $\int_F |U_n|^2\,d\mathcal A=1$ and $D^{\tau_n} U_n=s_n^2 U_n$ in $\HH$. 
	
	Under these conditions, the semiclassical measure of sequence $\{(U_n,1/\tau_n)\}_{n=1}^\infty$ is scalar multiple of $(\phi_1^{-1})_\sharp\mu_L$, where $\mu_L$ is the uniform Liouville measure on $S_1\HH$.
\end{predl}

\noindent (To define Liouville measure on $S_1\HH$, introduce coordinates therein: if vector $v\in S_1\HH$ has basepoint $z=x+iy\in\mathbb C^+$ and has oriented angle $\theta$ with geodesic line $\Ree z=\const$ then we take $(x,y,\theta)$ as coordinates of $v$. In these coordinates, $d\mu_L=\dfrac{dx\,dy\,d\theta}{2\pi y^2}$.)

\noindent {\bf Remark.} Proposition \ref{prop:horocyclic_que} surely applies to functions 
$$
U_{s,\tau} = \Exc^{(s)}_{0\to \tk} u_0^{(s)} = \frac{K_{\tau-1}}{\sqrt{s^2+\tau \cdot(\tau-1)}}\cdot\,\cdots\,\cdot\frac{K_1}{\sqrt{s^2+1\cdot(1+1)}}\cdot\frac{K_0}{\sqrt{s^2+0\cdot(0+1)}}\, u^{(s)}_0,
$$
if $\tau=[\tk s]$ grows faster than $s$, that is, if $\tk\to\infty$. Thus, Proposition \ref{prop:horocyclic_que} implies Theorem \ref{th:horrorcyclic_chaos}.

\medskip

\noindent {\bf Remark.}
We do not assume that $s_n\xrightarrow{n\to\infty}\infty$. This is because now we quantize at level $1/\tau_n$, this should necessarily go to $0$ with $n\to\infty$  since 	${\tau_n}/{s_n}\xrightarrow{n\to\infty}\infty$ and $s_n$ are separated from zero. (The latter is because spectra of all operators $-\Delta_{\HH}+2i\tau y\frac{\partial}{\partial x}$ ($\tau\in\mathbb Z$)  on $\mathcal F^{\tau}(\Gamma)$ differ one from another only by finite number of points from $\mathbb Z/4$, see \cite{Fay}.)

\medskip

\noindent {\bf Proof of Proposition \ref{prop:horocyclic_que}.} Let $\bar{\mu}^\infty$ be a semiclassical measure of a subsequence in $\{(U_n,1/\tau_n)\}_{n=1}^\infty$. Since $\dfrac{1}{\tau^2}\cdot D^\tau = -\dfrac1{\tau^2}\cdot\Delta_{\HH}+2i y\cdot\dfrac{1}{\tau}\cdot\dfrac{\partial}{\partial x} = \Op^{\R^2}_{1/\tau}(2H_1-1)$ and $\dfrac{1}{\tau_n^2}\cdot D^{\tau_n} U_n=c_n\cdot U_n$ with real $c_n=(s_n/\tau_n)^2\xrightarrow{n\to\infty}0$, measure $\bar{\mu}^\infty$ is concentrated on the set $\{H_1=1/2\}\subset T^*\HH$, cf. first assertion of Lemma \ref{lemma:Wigner_usual}. Arguing as in the proof of the second assertion of Lemma \ref{lemma:Wigner_usual}, we conclude that $\bar{\mu}^\infty$ is invariant under restriction of flow $\exp t\Xi_1$ onto this level set.

Therefore, measure $\underline{\mu}^\infty:=(\phi_1)_\sharp\bar{\mu}^\infty$ is concentrated on $S_1\HH$ and is invariant under action of horocyclic flow  $h^\infty$, we have used second and fourth assertions of Proposition~\ref{predl:t_tstar_change}. Applying Lemma \ref{lemma:invariance} for $\hbar_n=1/\tau_n$ and $\tk=1$, we conclude that $\underline{\mu}^\infty$ is $\Gamma$-invariant. By our assumptions, $X=\Gamma\setminus\HH$ is smooth compact surface. By Furstenberg's Theorem (\cite{Furstenberg73}, \cite{Marcus}), there exists only one Borel probability measure on $S_1X$ invariant under horocyclic flow, up to multiplicative constant, this property is known as \emph{unique ergodicity} of $h^\infty$. The proof is complete. $\blacksquare$

\medskip

\noindent To conclude with, we get rid of cotangent bundle and derive the following 

\begin{sled}
	Under assumptions of Proposition \ref{prop:horocyclic_que}, measures $|U_n|^2\cdot\mathcal A$ converge weak* to $\dfrac{\mathcal A}{\mathcal A(F)}$.
\end{sled}

\bigskip

\pagebreak

\noindent {\bf Acknowledgments.} Author is grateful to:
\begin{itemize}[label={---}]
	\item P.G. Zograf for introducing to the mysterious world of spectral hyperbolic geometry;
	\item A. Logunov for the idea to separate variables on cylinders and also for conversations leaded author to studying the horocyclic flow;
	\item R. Romanov for encouraging to investigate singularity propagation for waves traveling without Hamiltonian;
	\item all the \emph{collective} of Chebyshev Laboratory for creative atmosphere;
	\item the reviewers whose attention allowed to make this paper better.
\end{itemize}

\end{document}